\documentclass[leqno,12pt]{article}
\usepackage{amsmath,amssymb,rotating,a4wide,color}
\usepackage{wasysym}
\usepackage{hyperref}
\usepackage{adjustbox}

\usepackage[all,cmtip]{xy}

\usepackage{verbatim}

\newdir^{ (}{{}*!/-3pt/\dir^{(}}
\newdir_{ (}{{}*!/-3pt/\dir_{(}}

\entrymodifiers={+!!<0pt,\fontdimen22\textfont2>}

\def\Agemotext{\rotatebox[origin=c]{180}{$\textstyle\kern-1pt\varOmega$}}
\def\Agemoscript{\rotatebox[origin=c]{180}{$\scriptstyle\kern-1pt\varOmega\kern1pt$}}
\def\Agemoscriptscript{\rotatebox[origin=c]{180}{$\scriptscriptstyle\kern-1pt\varOmega\kern1pt$}}

\newcommand{\rtuple}{{r}}
\def\CAN{\Theta}
\let\FEpi\varphi
\let\XrRedMap\Phi
\def\strong{{\rm strong}}

\def\Sch{\mathop{\rm Sch}\nolimits}

\def\bigtimesdisplay{\mathop{\raise-2pt\hbox{\huge$\times$}}}
\def\bigtimestext{\mathop{\raise-1pt\hbox{\Large$\times$\kern-2pt}}}
\def\bigtimes{\mathchoice{\bigtimesdisplay}{\bigtimestext}{\bigtimestext}{\bigtimestext}}

\let\oldbigwedge\bigwedge
\def\newbigwedge{\mathord{\adjustbox{valign=B,totalheight=8.5pt}{$\oldbigwedge$}}}
\renewcommand{\bigwedge}{\newbigwedge}

\newbox\circbulletbox
\setbox\circbulletbox=\hbox{\,{\large$\circledcirc$}\kern-8.7pt\raise .6pt\hbox{$\bullet$}\,}

\let\le\leqslant
\let\ge\geqslant

\let\triangleleftnaked\triangleleft
\def\triangleleft{\mathrel{\triangleleftnaked}}
\def\depsilon{d^{\kern1pt\epsilon}}
\def\Fratt{\mathop{\rm Fr}\nolimits}
\def\circVbig{\hbox{\text{\it\r{V}}}}
\def\circVscript{\hbox{\scriptsize\text{\it\r{V}}}}
\def\circVscriptscript{\mbox{\tiny\text{\it\r{V}}}}

\def\circVlimits_#1^#2{{\mathchoice%
   {\circVbig{}^{\kern2pt #2}_{\kern-2pt #1}}%
   {\circVbig{}^{\kern2pt #2}_{\kern-2pt #1}}%
   {\scriptstyle\circVscript{}^{\kern1.7pt #2}_{\kern-1pt #1}}%
   {\scriptscriptstyle\circVscriptscript{}^{\kern1.5pt #2}_{\kern-1pt #1}}%
   }}
\def\circVr_#1{\circVlimits_#1^r}
\def\circVs_#1{\circVlimits_#1^s}

\def\circWbig{\hbox{\text{\it\r{W}}}}
\def\circWscript{\hbox{\scriptsize\text{\it\r{W}}}}
\def\circWscriptscript{\mbox{\tiny\text{\it\r{W}}}}

\def\circWlimits_#1^#2{{\mathchoice%
   {\circWbig{}^{\kern2pt #2}_{\kern-2pt #1}}%
   {\circWbig{}^{\kern2pt #2}_{\kern-2pt #1}}%
   {\scriptstyle\circWscript{}^{\kern1.7pt #2}_{\kern-1pt #1}}%
   {\scriptscriptstyle\circWscriptscript{}^{\kern1.5pt #2}_{\kern-1pt #1}}%
   }}

\def\OM{\mathchoice
  {\rlap{\kern3.2pt$\overline{\phantom{L}}$}M}
  {\rlap{\kern3.2pt$\overline{\phantom{L}}$}M}
  {\rlap{\kern2.4pt$\scriptstyle\overline{\phantom{L}}$}M}
  {\rlap{\kern1.8pt$\scriptscriptstyle\overline{\phantom{L}}$}M}}

\def\mycirc{{\kern1pt\circ\kern2pt}}

\def\Aut{\mathop{\rm Aut}\nolimits}

\def\Hom{\mathop{\rm Hom}\nolimits}

\def\Gal{\mathop{\rm Gal}\nolimits}

\def\Stab{\mathop{\rm Stab}\nolimits}

\def\trace{\mathop{\rm tr}\nolimits}

\def\rank{\mathop{\rm rank}\nolimits}

\def\GL{\mathop{\rm GL}\nolimits}
\def\SL{\mathop{\rm SL}\nolimits}

\def\ab{{\rm ab}}

\let\phi\varphi

\def\theta{\vartheta}
\let\epsilon\varepsilon
\let\setminus\smallsetminus

\newtheorem{Thm}{Theorem}[section]
\newtheorem{Prop}[Thm]{Proposition}
\newtheorem{Lem}[Thm]{Lemma}

\newtheorem{Cor}[Thm]{Corollary}

\newtheorem{Def}[Thm]{Definition}

\newtheorem{Rem}[Thm]{Remark}
\newtheorem{Ex}[Thm]{Example}

\newtheorem{Cons}[Thm]{Construction}

\numberwithin{Thm}{section}
\def\UseTheoremCounterForNextEquation{\setcounter{equation}{\value{Thm}}\addtocounter{Thm}{1}}

\def\qed{{\hskip0pt\unskip\unskip\nobreak\hfil\penalty50
          \hskip1em\hbox{}\nobreak\hfil
           {$\square$}
          \parfillskip=0pt\finalhyphendemerits=0
          \par}\medskip}

\newenvironment{Proof}
               {\noindent{\bf Proof.}\ }
               {\qed}

\newenvironment{Proofof}[1]
               {\noindent{\bf Proof of #1.}\ }
               {\qed}

\newcommand{\BF}{{\mathbb{F}}}

\newcommand{\BQ}{{\mathbb{Q}}}

\newcommand{\BZ}{{\mathbb{Z}}}

\newcommand{\CG}{{\cal G}}
\newcommand{\CH}{{\cal H}}

\newcommand{\CO}{{\cal O}}

\newcommand{\CR}{{\cal R}}

\newbox\mybox
\def\arrover#1{\mathrel{
       \setbox\mybox=\hbox spread 1.4em
              {\hfil$\scriptstyle#1$\hfil}
       \vbox{\offinterlineskip\copy\mybox
             \hbox to\wd\mybox{\rightarrowfill}}}}

\def\larrover#1{\mathrel{
       \setbox\mybox=\hbox spread 1.4em
              {\hfil$\scriptstyle#1\vphantom{g}$\hfil}
       \vbox{\offinterlineskip\copy\mybox
             \hbox to\wd\mybox{\leftarrowfill}}}}

\def\ontoover#1{\mathrel{
       \setbox\mybox=\hbox spread 1.4em
              {\hfil$\scriptstyle#1\vphantom{g}$\hfil}
       \vbox{\offinterlineskip\copy\mybox
             \hbox to\wd\mybox{\rightarrowfill\hskip-2.8mm
                               $\rightarrow$}}}}
\def\leftontoover#1{\mathrel{
       \setbox\mybox=\hbox spread 1.4em
              {\hfil$\scriptstyle#1\vphantom{g}$\hfil}
       \vbox{\offinterlineskip\copy\mybox
             \hbox to\wd\mybox{$\leftarrow$\hskip-2.8mm
                               \leftarrowfill}}}}
\let\longto\longrightarrow
\let\into\hookrightarrow
\let\onto\twoheadrightarrow
\def\longonto{\ontoover{\ }}

\def\isoto{\mathrel{
       \setbox\mybox=\hbox spread 0.9em
              {\hfil$\scriptstyle\sim$\hfil}
       \vbox{\offinterlineskip\copy\mybox
             \hbox to\wd\mybox{\rightarrowfill}}}}

\def\invlim{\mathop{\vtop{\hbox{\rm lim}\vskip-8pt
        \hbox{\hskip1pt$\scriptstyle\longleftarrow$}\vskip-1pt}}}

\begin{document}

\title{A Cohen-Lenstra Heuristic for Schur $\sigma$-Groups}

\author{\qquad
\begin{minipage}{.3\hsize}
Richard Pink\\[12pt]
\small Department of Mathematics \\
ETH Z\"urich\\
8092 Z\"urich\\
Switzerland \\
pink@math.ethz.ch\\[9pt]
\end{minipage}
\qquad
\begin{minipage}{.3\hsize}
Luca \'Angel Rubio\\[12pt]
\small Department of Mathematics \\
ETH Z\"urich\\
8092 Z\"urich\\
Switzerland \\
lrubio@ethz.ch\\[9pt]
\end{minipage}
\qquad
}

\date{\today}

\maketitle
\begin{center}
In memory of Nigel Boston and Giorgio Frangoni
\end{center}
\medskip
\begin{abstract}
For any odd prime $p$ and any imaginary quadratic field~$K$, the $p$-tower group $G_K$  associated to~$K$ is the Galois group over $K$ of the maximal unramified pro-$p$-extension of~$K$. This group comes with an action of a finite group $\{1,\sigma\}$  of order~$2$ induced by complex conjugation and is known to possess a number of other properties, making it a so-called Schur $\sigma$-group. Its maximal abelian quotient is naturally isomorphic to the $p$-primary part of the narrow ideal class group of~$\CO_K$, and the Cohen-Lenstra heuristic gives a probabilistic explanation for how often this group is isomorphic to a given finite abelian $p$-group.

The present paper develops an analogue of this heuristic for the full group~$G_K$.
It is based on a detailed analysis of general pro-$p$-groups with an action of $\{1,\sigma\}$, which we call $\sigma$-pro-$p$-groups. We construct a probability space whose underlying set consists of $\sigma$-isomorphism classes of weak Schur $\sigma$-groups and whose measure is constructed from the principle that the relations defining $G_K$ should be randomly distributed according to the Haar measure. We also compute the measures of certain basic subsets, the result being inversely proportional to the order of the $\sigma$-automorphism group of a certain finite $\sigma$-$p$-group, as has often been observed before. Finally, we show that the $\sigma$-isomorphism classes of weak Schur $\sigma$-groups for which each open subgroup has finite abelianization form a subset of measure~$1$.
\end{abstract}

{\renewcommand{\thefootnote}{}
%
}

\newpage
\renewcommand{\baselinestretch}{0.6}\normalsize
\tableofcontents
\renewcommand{\baselinestretch}{1.0}\normalsize

\section{Introduction}
\label{Intro}

Consider an odd prime $p$ and an imaginary quadratic field~$K$, and let $K_p$ denote the maximal unramified pro-$p$-extension of~$K$. Then $K_p$ is Galois over~$\BQ$, and its Galois group is the semidirect product of the pro-$p$-group $G_K := \Gal(K_p/K)$ with a group $\{1,\sigma\}$  of order~$2$ generated by complex conjugation. The group $G_K$ with the action of $\sigma$ is called the \emph{$p$-tower group} associated to~$K$.

By abelian class field theory the maximal abelian quotient of $G_K$ is naturally isomorphic to the $p$-primary part of the narrow ideal class group of~$\CO_K$. Cohen and Lenstra \cite{CohenLenstra1983} described a probabilistic explanation for how often this group is isomorphic to a given finite abelian $p$-group and provided evidence for it. 
Since then there has been much work on extending the Cohen-Lenstra heuristic to other cases. The present paper develops an analogue of this heuristic for the full group~$G_K$. For this we follow the same principles as Boston-Bush-Hajir \cite{BostonBushHajir2017}, but with a few notable differences. A similar approach is used in a more general setting in the very recent article \cite{LiuWoodZureickBrown2024} by Liu-Wood-Zureick-Brown.

\medskip
On the one hand we found it helpful to independently develop the concept of \emph{$\sigma$-pro-$p$-groups}, which are pro-$p$-groups together with an action of a fixed group $\{1,\sigma\}$ of order~$2$. 
These are interesting animals in their own right: 
Define the \emph{even} part of a $\sigma$-pro-$p$-group~$G$ as $G^+ := \{ a\in G \mid {}^\sigma a = a\}$ and the \emph{odd} part as $G^- := \{ a\in G\mid {}^\sigma a = a^{-1}\}$. Then the multiplication map induces a homeomorphism $G^+\times G^- \to G$, which therefore constitutes a kind of $\BZ/2\BZ$-grading of~$G$.
Also, the formation of these parts is exact in short exact sequences. 
(These facts have already been used essentially by Koch-Venkov in \cite{KochVenkov1974}.)
Moreover, under the homeomorphism $G^+\times G^- \to G$ the Haar measure on $G$ corresponds to the product of the Haar measure on the subgroup $G^+$ and a natural measure on~$G^-$, which we call the Haar measure on~$G^-$. To obtain probability measures we normalize all these measures to have total measure~$1$.

Furthermore, we define \emph{$\sigma$-isomorphisms} of $\sigma$-pro-$p$-groups as a $\sigma$-equivariant topological isomorphisms. 
We define \emph{$\sigma$-automorphisms} accordingly and denote the group of $\sigma$-automorphisms of $G$ by $\Aut_\sigma(G)$. 
A special role in our investigations is played by the $\sigma$-pro-$p$-group $F_n$ whose underlying pro-$p$-group is free on $n<\infty$ odd generators. As ingredient for one of our main results 
we establish some representation theoretic properties of open subgroups of~$F_n$.
We also need some technical facts concerning the group $\Aut_\sigma(F_n)$.

\medskip
Returning to $p$-tower groups, our goal is then to formulate a plausible heuristic for the $\sigma$-isomorphism class of the $\sigma$-pro-$p$-group~$G_K$. The main difference to \cite{BostonBushHajir2017} is that we work throughout with the full group $G_K$ instead of its finite quotients, and that we actually construct a probability space. 
This construction is based on the known fact that $G_K$ is $\sigma$-isomorphic to the quotient of $F_n$ by the closed normal subgroup $N_\rtuple$ generated by all conjugates of $n$ elements $r_1,\ldots,r_n$ of $F_n^-$ for some integer~$n$. We agree with \cite{BostonBushHajir2017} and \cite{LiuWoodZureickBrown2024} that these relations should be randomly distributed according to the Haar measure on~$F_n^-$.

\medskip
For any $p$-tower group $G_K$ one also knows that the abelianization of every open subgroup is finite, because it is isomorphic to the narrow ideal class group of a number field. We therefore had to decide whether to base our probability space on the set of $\sigma$-isomorphism classes of $\sigma$-pro-$p$-groups that satisfy this extra condition or not. We ultimately found it more natural not to carry along this condition and to study its effect only at the end. 

Also, following Koch-Venkov \cite{KochVenkov1974} the $\sigma$-pro-$p$-groups of the form $F_n/N_\rtuple$ described above and whose abelianizations are finite are already called \emph{Schur $\sigma$-groups} (see Proposition \ref{SchurSigmaProp}). This caused a little terminological dilemma in that our construction naturally yields a probability space without that condition, while the $p$-tower groups that we are really interested in satisfy an even stronger condition. We therefore opted for calling all $\sigma$-pro-$p$-groups of the form $F_n/N_\rtuple$ \emph{weak Schur $\sigma$-groups} and all those for which every open subgroup has finite abelianization \emph{strong Schur $\sigma$-groups}. Thus every strong Schur $\sigma$-group is a Schur $\sigma$-group, and every Schur $\sigma$-group is a weak Schur $\sigma$-group, though not vice versa. Weak Schur $\sigma$-groups up to $\sigma$-isomorphism can also be characterized abstractly by Proposition~\ref{WSSGEquiv}.

\medskip
We thus define our probability space as the set $\Sch$ (pronounce \emph{Schur}) of all $\sigma$-iso\-mor\-phism classes $[G]$ of weak Schur $\sigma$-groups $G\cong F_n/N_\rtuple$ for all~$n$. For any $[G]\in\Sch$ and any $i\ge2$ we consider the subset
$$U_{i,G}\ :=\ \bigl\{ [H]\in\Sch \bigm| H/D_i(H) \cong G/D_i(G) \bigr\},$$
where $D_i(\ )$ denotes the $i$-th step in the Zassenhaus filtration of a pro-$p$-group. (Any other canonical filtration cofinal to it, such as the descending $p$-central filtration, would do the same job.) These subsets form a basis of a unique topology on~$\Sch$, which is Hausdorff and totally disconnected. 
For any $n\ge0$ let $\Sch_n$ denote the subset of all $[G]\in\Sch$ for which the minimal number of generators of $G$ is equal to~$n$. Then $\Sch_n$ is an open and closed compact subset of $\Sch$ and therefore a profinite topological space.

\medskip
Next, for every integer $n\ge0$ there is a natural map
$$\CAN_n\colon\ (F_n^-)^n\longto\Sch,\ \ \rtuple = (r_1,\ldots,r_n) \mapsto [F_n/N_\rtuple].$$
We show that this map is continuous and can therefore define a probability measure $\mu_n$ on $\Sch$ as the pushforward under $\CAN_n$ of the Haar measure on $(F_n^-)^n$. This measure is supported on the subset $\Sch_{\le n} := \bigcup_{m\le n}\Sch_m$ and depends on~$n$. But for any fixed $m$, the restrictions $\mu_n|\Sch_m$ for all $n\ge m$ agree up to explicit proportionality factors depending on $m$ and~$n$. Using this fact we can prove that the measures $\mu_n$ converge for $n\to\infty$ to a canonical probability measure $\mu_\infty$ on~$\Sch$, whose restriction to each $\Sch_m$ is an explicit multiple of $\mu_m|\Sch_m$ (see Theorem \ref{MuinfMum}).

\medskip
Having established our probability space, we next compute the measures of certain subsets such as $U_{i,G}$. There being no reason to favor $G/D_i(G)$ over any other canonical finite quotient of~$G$, we do this more generally as follows. Fix an integer $n\ge0$ and a $\sigma$-invariant open subgroup $D$ of~$F_n$ that is contained in the Frattini subgroup of $F_n$ and that is invariant under $\Aut_\sigma(F_n)$. Then any $\sigma$-isomorphism class $[G] = [F_n/N_\rtuple]\in\Sch_n$ determines a unique $\sigma$-isomorphism class of finite $\sigma$-$p$-groups $G_D \cong F_n/N_\rtuple D$. Generalizing the definition of $U_{i,G}$ we are interested in the set 
$$U_{D,G}\ :=\ \bigl\{ [H]\in\Sch \bigm| H_D \cong G_D \bigr\}.$$
In Proposition \ref{UDGMeasureOpen} we compute its measure as
$$\mu_\infty(U_{D,G})\ =\ \frac{C_\infty}{C_{n-m_{D,G}}}\cdot \frac{1}{\bigl|\Aut_\sigma(G_D)\bigr|},$$
where $m_{D,G}$ is the minimal number of relations defining $G_D$ as a quotient of $F_n/D$ and the constants are positive and given by explicit formulas. This agrees with the formula for $u=0$ from Liu-Wood-Zureick-Brown \cite[Thm.\,7.4]{LiuWoodZureickBrown2024}, and in the special case that $D$ is a step in the descending $p$-central filtration of~$F_n$ it agrees with the formula from Boston-Bush-Hajir \cite[Conj.\,1.3]{BostonBushHajir2017}, taking into account the proportionality factor from Proposition~\ref{MuinfOfSchn}.

In the case that $G$ itself is finite, in Proposition \ref{FinSchurMeas} we show that the singleton $\{[G]\}$ is open and closed in $\Sch$ of measure
$$\mu_\infty(\{[G]\})\ =\ \frac{C_\infty}{\bigl|\Aut_\sigma(G)\bigr|}.$$
Again this agrees with the formula from \cite[Conj.\,1.3]{BostonBushHajir2017}, taking Proposition \ref{MuinfOfSchn} into account. Propositions \ref{UDGMeasureFin} and \ref{UDGMeasureFinVar} provide corresponding results under the weaker assumption that $G_D$ is finite.

\medskip
As explained above, all these results concern $\sigma$-isomorphism classes of weak Schur $\sigma$-groups, whereas $p$-tower groups are strong Schur $\sigma$-groups. The heuristic can therefore only be plausible if the subset $\Sch^\strong$ of $\sigma$-isomorphism classes of strong Schur $\sigma$-groups has measure~$1$. 
(This question was also asked more generally in \cite[\S1.5]{LiuWoodZureickBrown2024}.)
That we can indeed prove this in Theorem \ref{StrongThm} would seem to indicate that our approach is on the right track.

\medskip
A further interesting question is under which conditions and how often a $p$-tower group is finite. After Golod and Shafarevich \cite{GolodShafarevich1964} proved that $G_K$ can be infinite, there has been much work on this, and most known results concern abstract Schur $\sigma$-groups. 

So let $G$ be a strong Schur $\sigma$-group with minimal number of generators~$d_G$. According to current knowledge $G$ is finite if $d_G\le1$, infinite if $d_G\ge3$, and in the case $d_G=2$ it is infinite unless its Zassenhaus type is (3,3) or (3,5) or (3,7) (see McLeman \cite[\S2]{McLeman2008}). In each of those remaining cases it is known that the group can be finite, but the cases are not yet settled completely. McLeman \cite[Conj.\,2.9]{McLeman2008}  conjectured that any $p$-tower group of Zassenhaus type (3,3) should be finite. 
Partial evidence for this is provided in a separate paper \cite[Prop.\,4.7]{Pink2025} by the first author, proving that the set of $[G]\in\Sch$ such that $G$ is infinite of Zassenhaus type (3,3) has measure~$0$ in the case $p>3$. 

\medskip
Observe that such a statement naturally requires a probability space based on possibly infinite (strong) Schur $\sigma$-groups and not only on some finite quotients thereof. This approach also allows us to formulate properties of finite quotients within one coherent theory. We therefore view our work as a natural next step from \cite{BostonBushHajir2017} and hope that it will bring other benefits in the future.


\medskip
The article is structured as follows: In Section \ref{PPG} we recall known facts about pro-$p$-groups. 
In the subsequent four sections we study $\sigma$-pro-$p$-groups from different angles: In Section \ref{SPG} we establish basic properties of finite $\sigma$-$p$-groups. In Section \ref{SPPG} we generalize these to $\sigma$-pro-$p$-groups and discuss the Haar measures on them. In Section \ref{FPPG} we study some representation theoretic aspects of open subgroups of $\sigma$-pro-$p$-groups whose underlying pro-$p$-groups are free on finitely many generators. 
(This is needed for the proof of Theorem \ref{StrongThm}.)
In Section \ref{AutFSG} we briefly discuss the group of $\sigma$-automorphisms of a $\sigma$-pro-$p$-group $F_n$ whose underlying pro-$p$-group is free on $n<\infty$ odd generators.

In the next three sections we construct our probability space. In Section \ref{WSSG} we introduce the central concept of weak Schur $\sigma$-groups and give equivalent characterizations. In Section \ref{SchTop} we construct the topology on $\Sch$ and derive its basic properties. In Section \ref{ProbMeas} we then construct the probability measure $\mu_\infty$ on~$\Sch$.

In Section \ref{CompProb} we derive explicit formulas for the measures of certain basic subsets of~$\Sch$. 
In Section \ref{SSSG} we prove that the subset $\Sch^\strong$ of $\sigma$-isomorphism classes of strong Schur $\sigma$-groups has measure~$1$. 
Finally, in Section \ref{pTG} we discuss the application to $p$-tower groups and 
to the question how often they are finite.

\medskip
This article grew out of the master thesis of the second author that was supervised by the first author.

\section{Pro-$p$-groups}
\label{PPG}

Throughout this article, all homomorphisms of pro-$p$-groups are tacitly assumed to be continuous, and all subgroups are tacitly assumed to be closed. Thus by the subgroup generated by a subset of a pro-$p$-group we always mean the closure of the abstract subgroup generated by that set. In particular, when we say that a pro-$p$-group is generated by certain elements, we mean that it is topologically generated by them. 

\medskip
Consider a pro-$p$-group~$G$. For any elements $a,b\in G$ we abbreviate ${}^ab := aba^{-1}$ and $[a,b] := aba^{-1}b^{-1}$, and we write $a\sim b$ if and only if the elements are conjugate. 
For any elements $a_1,\ldots a_n\in G$ we let $\langle a_1,\ldots,a_n\rangle$ denote the subgroup generated by $a_1,\ldots a_n$.
For any subsets $A,B\subset G$ we let $[A,B]$ denote the subgroup generated by the subset $\{[a,b]\mid a\in A,\; b\in B\}.$ 
The notation $H<G$ means that $H$ is a subgroup of~$G$, but not necessarily a proper subgroup.
Likewise the notation $A\subset B$ means that $A$ is a subset of~$B$, but not necessarily a proper subset.

\medskip
For any $i\ge1$ we let $D_i(G)$ denote the $i$-th dimension subgroup of $G$ (see \cite[\S11]{DdSMS2003}). These constitute the Zassenhaus filtration of~$G$ and are characteristic subgroups of $G$ whose intersection is the trivial subgroup $\{1\}$. In particular every open subgroup of $G$ contains some $D_i(G)$. Also, 
any homomorphism $G\to H$ induces a homomorphism $D_i(G)\to D_i(H)$, 
and for any normal subgroup $N\triangleleft G$ we have $D_i(G/N) = D_i(G)N/N$. If $G$ is finitely generated, then each $D_i(G)$ is an open subgroup of~$G$ and every open subgroup of $G$ is finitely generated.
%

\medskip
The subgroup $D_2(G)$ is equal to the Frattini subgroup $\Fratt(G)$ of~$G$, so that $G/\Fratt(G)$ is the maximal quotient of $G$ that is an $\BF_p$-vector space. For later use we record the following direct consequences of \cite[Prop.\,3.9.1, Cor.\,3.9.3]{NSW2008}:

\begin{Prop}\label{NormGens}
\begin{enumerate}
\item[(a)] For any subgroup $H<G$ with $H\Fratt(G)=G$ we have $H=G$.
\item[(b)] For any normal subgroups $M\triangleleft G$ and $N\triangleleft G$ with $M\Fratt(N)[G,N]=N$ we have $M=N$.
\end{enumerate}
\end{Prop}

Next, a pro-$p$-group $G$ is called free on finitely many elements $x_1,\ldots,x_n \in G$ if it satisfies the usual universal property for homomorphisms into all pro-$p$-groups~$H$. This means that for any choice of $h_1,\ldots,h_n\in H$ there exists a unique homomorphism $\phi\colon {G\to H}$ with $\phi(x_i)=h_i$ for all~$i$. The residue classes of $x_1,\ldots,x_n$ then form an $\BF_p$-basis of $G/\Fratt(G)$. In particular $x_1,\ldots,x_n\in G$ are distinct and generate $G$ by Proposition \ref{NormGens} (a). Also the number $n$ is uniquely determined by~$G$ up to isomorphism and called the rank of~$G$.

\begin{Prop}\label{AutFGCons}
If $G$ is free on $x_1,\ldots,x_n \in G$, then for any $y_1,\ldots,y_n\in G$ whose residue classes form an $\BF_p$-basis of $G/\Fratt(G)$, there exists a unique automorphism $\alpha$ of $G$ with $\alpha(x_i)=y_i$ for all~$i$. In particular $G$ is free on the elements $y_1,\ldots,y_n$.
\end{Prop}

\begin{Proof}
By the universal property there exists a unique endomorphism $\alpha\colon G\to G$ with $\alpha(x_i)=y_i$ for all~$i$. Proposition \ref{NormGens} (a) implies that $\alpha$ is surjective. For any $i\ge0$ it thus induces a surjective endomorphism of the finite group $G/D_i(G)$, which is therefore an automorphism, and so $\alpha$ is an automorphism of $G = \invlim_i G/D_i(G)$. This in turn implies that the universal property also holds for the elements $y_1,\ldots,y_n$.
\end{Proof}

\medskip
Finally, there is a natural isomorphism 
\UseTheoremCounterForNextEquation
\begin{equation}\label{H1Cong}
H^1(G,\BF_p)\ \cong\ \Hom(G/\Fratt(G),\BF_p).
\end{equation}
The dimension $d_G := \dim_{\BF_p}\!H^1(G,\BF_p)$ is the minimal number of generators of~$G$ (see \cite[Prop.\,3.9.1]{NSW2008}).
Thus if $F_n$ denotes the free pro-$p$-group of rank~$n$, we have $d_G\le n$ if and only if $G\cong F_n/N$ for some normal subgroup $N\triangleleft F_n$. In that case we have $d_G=n$ if and only if $N\subset\Fratt(F_n)$.
Assuming that, there is a natural isomorphism
\UseTheoremCounterForNextEquation
\begin{equation}\label{H2Cong}
H^2(G,\BF_p)\ \cong\ H^1(N,\BF_p)^{F_2}\ \cong\ \Hom\bigl(N/\Fratt(N)[G,N],\BF_p\bigr),
\end{equation}
and $\dim_{\BF_p}\!H^2(G,\BF_p)$ is the minimal number of generators of~$N$ as a normal subgroup of $F_n$ and hence the minimal number of relations of $G$ as a quotient of~$F_n$ (see \cite[Prop.\,3.9.5]{NSW2008}).

\section{$\sigma$-$p$-Groups}
\label{SPG}

Throughout this article we fix a prime $p>2$ and a finite group $\{1,\sigma\}$ of order~$2$. We call a finite $p$-group $G$ with an action of $\{1,\sigma\}$ a \emph{$\sigma$-$p$-group} and denote the action by $a\mapsto{}^\sigma a$. A $\sigma$-equivariant homomorphism between $\sigma$-$p$-groups is called a \emph{$\sigma$-homomorphism}. This defines a \emph{category of $\sigma$-$p$-groups} which possesses kernels and images and short exact sequences. 
The notions \emph{$\sigma$-isomorphisms} and \emph{$\sigma$-automorphisms} are defined in the evident way. The group of $\sigma$-automorphisms of a $\sigma$-$p$-group $G$ is denoted $\Aut_\sigma(G)$.

\medskip
Consider a $\sigma$-$p$-group~$G$. Then for any $\sigma$-invariant normal subgroup $N\triangleleft G$ the factor group $G/N$ is a $\sigma$-$p$-group and the projection map $G\onto G/N$ is a $\sigma$-homomorphism. In particular this holds for every characteristic subgroup of $G$, for instance for the center $Z(G)$ and the commutator subgroup $[G,G]$ with the abelianization $G_\ab := G/[G,G]$. 

\medskip
To $G$ we associate the natural subsets
\UseTheoremCounterForNextEquation
\begin{eqnarray}\label{G+Def}
G^+\ :=\ G^{+1} &\!:=\!& \bigl\{ a\in G\bigm| {}^\sigma a = a\bigr\}, \\[3pt]
\UseTheoremCounterForNextEquation\label{G-Def}
G^-\ :=\ G^{-1} &\!:=\!& \bigl\{ a\in G\bigm| {}^\sigma a = a^{-1}\bigr\}.
\end{eqnarray}
We call the elements of $G^+$ \emph{even} and the elements of $G^-$ \emph{odd}. We call $G$ \emph{totally even} if $G=G^+$, and \emph{totally odd} if $G=G^-$.

\medskip
While $G^+$ is always a subgroup of~$G$, the subset $G^-$ is in general not, though it is invariant under conjugation by~$G^+$ by direct computation. Also for any $g\in G$ one quickly shows that $g\cdot{}^\sigma g^{-1} \in G^-$. Both subsets are invariant under~$\sigma$, and any $\sigma$-homomorphism $G\to H$ induces a homomorphism $G^+\to H^+$ and a map $G^-\to H^-$. 

\medskip
If $G$ is abelian, for each $\epsilon\in\{\pm1\}$ the subset $G^\epsilon$ is simply the eigenspace for the eigenvalue $\epsilon$ of~$\sigma$ and thus a subgroup. Since $|G|$ is odd, the inclusions then induce a natural isomorphism $G^+\times G^-\isoto G$ that is functorial in~$G$. An abelian $\sigma$-$p$-group is thus the same as a $\BZ/2\BZ$-graded abelian $p$-group with even part $G^+$ and odd part~$G^-$.

\medskip
For induction proofs the following facts are crucial:

\begin{Prop}\label{SPGIndLem}
Every non-trivial $\sigma$-invariant normal subgroup $N\triangleleft G$ possesses a $\sigma$-invariant subgroup of order~$p$ that is normal in~$G$, and a $\sigma$-invariant subgroup of index~$p$ that is normal in~$G$.
\end{Prop}

\begin{Proof}
Since $N$ is a non-trivial normal subgroup of the $p$-group~$G$, the subgroup $N\cap Z(G)$ is non-trivial. As this is an abelian $\sigma$-$p$-group, at least one of its even or odd parts is non-trivial. Any element of order $p$ thereof thus generates a $\sigma$-invariant subgroup of $N$ of order~$p$ that is normal in~$G$.

Similarly, since $N$ is a non-trivial normal subgroup of the $p$-group~$G$, the subgroup $[G,N]$ is a proper subgroup of~$N$. The factor group $N/[G,N]$ is then a non-trivial abelian $\sigma$-$p$-group, and so at least one of its even or odd parts is non-trivial. Choosing a subgroup of index $p$ thereof, its inverse image under the projection from $N$ to the respective part of $N/[G,N]$ is a $\sigma$-invariant subgroup of index~$p$ that is normal in~$G$.
\end{Proof}


\begin{Prop}\label{SPGExact}
For any $\sigma$-invariant normal subgroup $N\triangleleft G$ and any $\epsilon\in\{\pm1\}$ the projection $G\onto G/N$ induces a surjective map $G^\epsilon \onto (G/N)^\epsilon$ all of whose fibers have the same cardinality $|N^\epsilon|$. In particular we have
$$|G^\epsilon|\ =\ |(G/N)^\epsilon|\cdot |N^\epsilon|.$$
\end{Prop}

\begin{Proof}
If $G$ is abelian, the first statement follows from the fact that the functor $G\mapsto G^\epsilon$ on abelian $\sigma$-$p$-groups is exact. For the general case we use induction on $|N|$, the case $|N|=1$ being trivial. 
If $|N|=p$, then $N$ is contained in the center of~$G$. For every element $\bar g\in (G/N)^\epsilon$, the inverse image of the subgroup $\langle\bar g\rangle$ is then an abelian $\sigma$-$p$-group. By the abelian case, the fiber of $\bar g$ thus has cardinality $|N^\epsilon|$ and is in particular non-empty, so the first statement follows for~$G$. 
If $|N|>p$, using Proposition \ref{SPGIndLem} we choose a $\sigma$-invariant subgroup $M<N$ of order~$p$ that is normal in~$G$. Then the statement is already proved for $M\triangleleft G$, and by the induction hypothesis it holds for $N/M\triangleleft G/M$. Together it then follows for $N\triangleleft G$, and we are done.

Finally, the cardinality formula is a direct consequence of the first statement.
\end{Proof}


\medskip
In the non-abelian case we still recover a kind of $\BZ/2\BZ$-grading:

\begin{Prop}\label{SPGGrad}
The product map induces bijections
$$\begin{array}{c}
G^+\times G^- \longto G,\ (a,b)\mapsto ab. \\[3pt]
G^-\times G^+ \longto G,\ (a,b)\mapsto ab.
\end{array}$$
In particular we have 
$$|G|\ =\ |G^+|\cdot|G^-|.$$
\end{Prop}

\begin{Proof}
We already know the first statement when $G$ is abelian, and this implies the cardinality formula in that case. In particular that formula holds whenever $|G|\le p$. Using the formula for $|G^\epsilon|$ from Proposition \ref{SPGExact} and induction on~$|G|$, the cardinality formula thus follows for all~$G$. By equality of cardinalities, to prove that the maps are bijective it remains to prove that they are surjective.

To achieve this we use induction on $|G|$. We already know the statement for $|G|\le p$, so assume that $|G|>p$. Using Proposition \ref{SPGIndLem} choose a $\sigma$-invariant normal subgroup $M\triangleleft G$ of order~$p$. Then by the induction hypothesis the maps are surjective for $G/M$ in place of~$G$. Every element $g\in G$ can therefore be written in the form $g=ab$ with $aM\in(G/M)^+$ and $bM\in(G/M)^-$. Using Proposition \ref{SPGExact} we can write $a=a'm$ and $b=nb'$ with $a'\in G^+$ and $b'\in G^-$ and $m,n\in M$. Thus $g = (a'mn)b' = a'(mnb')$.

Since $M$ has order~$p$, it is equal to $M^+$ or $M^-$. If $M=M^+$, the fact that $G^+$ is a subgroup implies that $a'mn\in G^+$, and the decomposition $g = (a'mn)b'$ does the job. If $M=M^-$, the fact that $M$ is contained in the center of $G$ implies that ${}^\sigma(mnb') = {}^\sigma(b'nm) = (b')^{-1}n^{-1}m^{-1} = (mnb')^{-1}$ and hence $mnb'\in G^-$. Thus in that case the decomposition $g = a'(mnb')$ does the job. By induction this shows that the first map is surjective. The analogous argument shows that the second map is surjective.
\end{Proof}

\medskip
We also have the following curious facts about conjugation. As we will not really use the first and its proof is analogous to that of the second, we leave the proof as an exercise.

\begin{Prop}\label{SPGConj}
For any $\epsilon\in\{\pm1\}$ and any $a\in G$ with ${}^\sigma a \sim a^\epsilon$ there exists $b\in G^\epsilon$ with $a\sim b$, and this $b$ is unique up to conjugation by~$G^+$.
\end{Prop}

\begin{Prop}\label{SPGSemidirect}
For any $a\in G$ there exists $b\in G^+$ such that $a\sigma$ is conjugate to $b\sigma$ in the semidirect product $G\rtimes\{1,\sigma\}$, and this $b$ is unique up to conjugation by~$G^+$.
\end{Prop}


\begin{Proof}
We prove this by induction on $|G|$, the case $|G|=1$ being trivial. Otherwise, using Proposition \ref{SPGIndLem} choose a $\sigma$-invariant normal subgroup $M\triangleleft G$ of order~$p$. Then by the induction hypothesis the proposition holds for $G/M$. After replacing $a\sigma$ by a conjugate we may therefore assume that the residue class of $a$ lies in $(G/M)^+$. 

Since $M$ has order~$p$, it is equal to $M^+$ or $M^-$. If $M=M^+$, then Proposition \ref{SPGExact} directly implies that $a\in G^+$, and $b:=a$ does the job. So assume that $M=M^-$. Since $M$ is contained in the center of~$G$, for every $m\in M$ we then have ${}^m(a\sigma) = ma\sigma m^{-1} = mam\sigma = am^2\sigma$. As $p$ is odd, the map $M\mapsto M$, $m\mapsto m^2$ is surjective, so it follows that $a\sigma$ is conjugate to $b\sigma$ for every element $b$ of the coset $aM$. Since one of these lies in $G^+$ by Proposition \ref{SPGExact}, we have proved the existence part in~$G$.

To prove the uniqueness part by induction, consider elements $a,b\in G^+$ such that $a\sigma$ is conjugate to $b\sigma$ under $G\rtimes\{1,\sigma\}$. Since $a\sigma$ lies in the non-trivial coset of~$G$, they are then already conjugate under~$G$. By the induction hypothesis the images of $a$ and $b$ in $G/M$ are conjugate under $(G/M)^+$. By Proposition \ref{SPGExact} we may thus choose an element $c\in G^+$ such that ${}^ca\equiv b$ modulo~$M$. Then $c\in G^+$ implies that ${}^c(a\sigma) = ({}^ca)\sigma$; so we may replace $a$ by ${}^ca$ and hence assume that $b=am$ for some $m\in M$. 

If $m=1$, we already have $a=b$ and are done. Otherwise, since $M$ has order~$p$, it is equal to $M^+$ or $M^-$ and generated by~$m$. If $M=M^-$, the fact that $a,b\in G^+$ and Proposition \ref{SPGExact} directly imply that $a=b$ again. So assume that $M=M^+=\langle m\rangle$. Consider the subgroup $H := \{d\in G\mid {}^d(a\sigma)\equiv a\sigma \bmod M\}$. Since $a\in G^+$ and hence ${}^\sigma(a\sigma)=a\sigma$, direct computation shows that $H$ is $\sigma$-invariant. Using the facts that $M^+=M\subset Z(G)$, two more direct computations show that we have a $\sigma$-equivariant homomorphism
$$\phi\colon H \longto M,\ \ d\mapsto {}^d(a\sigma)(a\sigma)^{-1}.$$
Since $M=M^+$, Proposition \ref{SPGExact} implies that $\phi(H)=\phi(H)^+=\phi(H^+)$. On the other hand, by assumption there exists an element $d\in G$ with ${}^d(a\sigma) = b\sigma = am\sigma$. Thus $d$ lies in $H$ and satisfies $\phi(d) = m$. Since $\phi(H)=\phi(H^+)$, there then also exists an element $d'\in H^+$ with $\phi(d') = m$. This element then satisfies ${}^{d'\kern-1pt}a{\cdot}\sigma = {}^{d'\kern-1pt}(a\sigma) = am\sigma = b\sigma$ and hence ${}^{d'\kern-1pt}a = b$. Thus $a$ and $b$ are conjugate under $G^+$, finishing the proof of the uniqueness part.
\end{Proof}

\section{$\sigma$-Pro-$p$-groups}
\label{SPPG}

We call a pro-$p$-group $G$ with an action of $\{1,\sigma\}$ a \emph{$\sigma$-pro-$p$-group}. All notions from Section \ref{SPG} have their counterpart here with the same notation and the same definitions. 

\medskip
Consider a $\sigma$-pro-$p$-group~$G$. Then $G^+$ is a closed subgroup and $G^-$ is a closed subset of~$G$, both endowed with the induced topology. Whenever $G$ is the filtered inverse limit of finite $\sigma$-$p$-groups~$G_i$, direct computation shows that $G^\epsilon$ is the corresponding inverse limit of the finite sets $G_i^\epsilon$ for every $\epsilon\in\{\pm1\}$.
Using this, Propositions \ref{SPGExact} and \ref{SPGGrad} imply:

\begin{Prop}\label{SPPGExact}
For any $\sigma$-invariant normal subgroup $N\triangleleft G$ and any $\epsilon\in\{\pm1\}$ the projection $G\onto G/N$ induces a surjective map $G^\epsilon \onto (G/N)^\epsilon$.
\end{Prop}

\begin{Prop}\label{SPPGGrad}
The product map induces homeomorphisms
$$\begin{array}{c}
G^+\times G^- \longto G,\ (a,b)\mapsto ab. \\[3pt]
G^-\times G^+ \longto G,\ (a,b)\mapsto ab.
\end{array}$$
\end{Prop}

Also, the results about conjugation from Section \ref{SPG} remain true in the pro-$p$ case:

\begin{Prop}\label{SPPGConj}
For any $\epsilon\in\{\pm1\}$ and any $a\in G$ with ${}^\sigma a \sim a^\epsilon$ there exists $b\in G^\epsilon$ with $a\sim b$, and this $b$ is unique up to conjugation by~$G^+$.
\end{Prop}

\begin{Proof}
Same as the proof of Proposition \ref{SPPGSemidirect} below with elements $c_i\in G_i$ such that ${}^{c_i}a_i\in G_i^\epsilon$ instead.
\end{Proof}

\begin{Prop}\label{SPPGSemidirect}
For any $a\in G$ there exists $b\in G^+$ such that $a\sigma$ is conjugate to $b\sigma$ in the semidirect product $G\rtimes\{1,\sigma\}$, and this $b$ is unique up to conjugation by~$G^+$.
\end{Prop}

\begin{Proof}
Write $G$ as the inverse limit of finite $\sigma$-$p$-groups~$G_i$ for $i$ in some filtered directed set~$I$. For every $i$ let $a_i\in G_i$ denote the image of~$a$. Then by Proposition \ref{SPGSemidirect} the set $S_i$ of $c_i\in G_i$ such that ${}^{c_i}(a_i\sigma)\in G_i^+\sigma$ is non-empty. Also, every morphism $\phi_{ji}\colon G_j\to G_i$ in the system induces a map $S_j\to S_i$. Since every filtered inverse limit of non-empty finite sets is non-empty, 
there exists an element $c\in G$ whose image in each $G_i$ lies in~$S_i$. This means that the image of ${}^c(a\sigma)$ in each $G_i\sigma$ lies in $G_i^+\sigma$, and so ${}^c(a\sigma)$ lies in $G^+\sigma$, proving the existence part. The uniqueness part is proved in the same way.
\end{Proof}

\medskip
Now recall that any pro-$p$-group $G$ possesses a unique Haar measure $\mu_G$ of total measure~$1$. This is the unique countably additive measure on the Borel $\sigma$-algebra of $G$ such that $\mu(gH) = [G:H]^{-1}$ for every $g\in G$ and every open subgroup $H<G$. Its construction, for instance in Fried-Jarden \cite[\S21]{FriedJarden2023}, does not actually require the group structure. All it needs is a filtered inverse limit $X$ of finite sets $X_i$ such that all fibers of the transition maps $X_j\to X_i$ have the same positive cardinality. By Proposition \ref{SPGExact} this is true for the system of sets $(G/N)^\epsilon$ for all $\sigma$-invariant open normal subgroups $N\triangleleft G$. In particular this defines a unique Borel measure $\mu_{G^\epsilon}$ of total measure $1$ that we call the \emph{Haar measure on~$G^\epsilon$}. For $\epsilon=+1$ this is of course just the usual Haar measure on the profinite group~$G^+$.

Moreover, Proposition \ref{SPGGrad} and the evident analogue of \cite[\S21.4]{FriedJarden2023} imply that the product measure $\mu_{G^+}\otimes\mu_{G^-}$ on $G^+\times G^-$ corresponds to $\mu_G$ under the homeomorphism from Proposition \ref{SPPGGrad}. Thus the measure $\mu_{G^-}$ can be constructed equivalently as the image measure of $\mu_G$ under the projection to the second factor.


\begin{Prop}\label{SPGQuot}
For any surjective $\sigma$-homomorphism of $\sigma$-pro-$p$-groups $\pi\colon G\onto H$  we have $\pi_*\mu_G = \mu_H$ and for any $\epsilon\in\{\pm1\}$ the induced map $\pi^\epsilon\colon G^\epsilon \to H^\epsilon$ satisfies $\pi^\epsilon_*\mu_{G^\epsilon} = \nobreak \mu_{H^\epsilon}$.
\end{Prop}

\begin{Proof}
Direct consequence of the construction of the respective Haar measures.
\end{Proof}

\begin{Prop}\label{SPGSubgroup}
Consider a $\sigma$-invariant open normal subgroup $N\triangleleft G$ and any $\epsilon\in\{\pm1\}$.
\begin{enumerate}
\item[(a)] The projection $G\onto G/N$ induces a surjective map $G^\epsilon \onto (G/N)^\epsilon$ all of whose fibers have the same measure $\mu_{G^\epsilon}(N^\epsilon)>0$.
\item[(b)] The restriction of $\mu_G$ to $N$ is equal to $|G/N|^{-1}\cdot \mu_N$.
\item[(c)] The restriction of $\mu_{G^\epsilon}$ to $N^\epsilon$ is equal to $|(G/N)^\epsilon|^{-1}\cdot \mu_{N^\epsilon}$.
\end{enumerate}
\end{Prop}

\begin{Proof}
All statements follow by direct computation from the definition of $\mu_{G^\epsilon}$ and Proposition \ref{SPGExact}.
\end{Proof}

\section{Free $\sigma$-pro-$p$-groups}
\label{FPPG}

In this section we study some representation theoretic aspects of open subgroups of free $\sigma$-pro-$p$-groups. We begin with free pro-$p$-groups without an action of~$\sigma$. First we recall Schreier's theorem in the case of pro-$p$-groups, see for instance \cite[Cor.\,3.9.6]{NSW2008}:

\begin{Prop}\label{FPPGOpenFree}
Any open subgroup $H$ of a finitely generated free pro-$p$-group $G$ is a finitely generated free pro-$p$-group and its minimal number of generators satisfies
$$d_H-1\ =\ [G:H]\cdot(d_G-1).$$
\end{Prop}

The following proposition makes this a little more precise in a special case:

\begin{Prop}\label{FPPGCyclicQuot}
Consider a finitely generated free pro-$p$-group $G$ on $n$ generators and an open normal subgroup $N$ with quotient $G/N \cong \BZ/p^r\BZ$ for $r\ge1$. Then there exist generators $x_1,\ldots,x_n$ of~$G$ such that $x_2,\ldots,x_n\in N$, and for any such choice $N$ is a free pro-$p$-group on the generators $x_1^{p^r}$ and $x_1^ix_jx_1^{-i}$ for all $0\le i<p^r$ and $2\le j\le n$.
\end{Prop}

\begin{Proof}
By assumption $G/\Fratt(G)$ is an $\BF_p$-vector space of dimension~$n$, and the image of $N$ is a subspace of dimension $n-1$. Choose a basis $\bar x_2,\ldots,\bar x_n$ of that subspace and an element $\bar x_1$ not in it. Lift these elements to elements $x_1\in G$ and $x_2,\ldots,x_n \in N$. Then Proposition \ref{AutFGCons} implies that $G$ is free on the generators $x_1,\ldots,x_n$.

Next let $\CG$ denote the abstract free group with generators $x_1,\ldots,x_n$ and set $\CH := \CG\cap H$. Then the subset $\CR := \{x_1^i\bigm| 0\le i<p^r\} \subset \CG$ is a set of representatives of the quotient $\CG/\CH$. For any $0\le i<p^r$ and $1\le j\le n$ the representative in $\CR$ of the coset $x_1^ix_jH$ is $x_1^i$ if $j\ge2$, respectively $x_1^{i+1}$ if $j=1$ and $i<p^r-1$, respectively $1$ if $j=1$ and $i=p^r-1$. By Schreier's lemma (see for instance \cite[Lemma 4.2.1]{Seress2003}) 
the subgroup $\CH$ is therefore generated by the elements $x_1^ix_jx_1^{-i}$ for all $0\le i<p^r$ and $2\le j\le n$ together with $1$ and $x_1^{p^r}$. 
Thus $\CH$ is generated by the elements in the proposition.

On the other hand, since $\CG$ is dense in~$G$, the fact that $H$ is open and closed in $G$ implies that $\CH$ is dense in~$H$. Thus $H$ is generated as a pro-$p$-group by the stated elements. Finally, this number of generators is $1+p^r(n-1)$. Since $H$ is free on this number of generators by Proposition \ref{FPPGOpenFree}, it is therefore free on the stated generators by Proposition \ref{AutFGCons}.
\end{Proof}

\begin{Thm}\label{FreeSubRep}
Let $G$ be a finitely generated free pro-$p$-group and let $N$ be an open normal subgroup of~$G$. Then for the action of $G/N$ on $N_\ab$ induced by the conjugation of $G$ on~$N$, there exists an isomorphism of $\BQ_p[G/N]$-modules
$$(N_\ab\otimes_{\BZ_p}\BQ_p) \oplus \BQ_p[G/N] \ \cong\ \BQ_p\oplus\BQ_p[G/N]^{\oplus d_G},$$
where $\BQ_p$ carries the trivial representation of~$G/N$.
\end{Thm}

\begin{Proof}
It suffices to show that the characters of the two representations coincide. So let $\chi$ denote the character of $N_\ab\otimes_{\BZ_p}\BQ_p$. Since $\BQ_p[G/N]$ is the regular representation of $G/N$, it suffices to show that 
\UseTheoremCounterForNextEquation
\begin{equation}\label{FreeSubRep1}
\chi([a])\ =\ {\scriptstyle\biggl\{}\!
\begin{array}{cl}
1 + [G:N]\cdot(d_G-1) & \hbox{if $a\in N$,} \\[3pt]
1 & \hbox{if $a\not\in N$,}
\end{array}
\end{equation}
for every $a\in G$. 

In the case $a\in N$ the residue class $[a]$ acts trivially on~$N_\ab$; hence $\chi([a])$ is the dimension of $N_\ab\otimes_{\BZ_p}\BQ_p$. Since $N$ is a free pro-$p$-group on $1+[G:N]\cdot(d_G-1)$ generators by Proposition \ref{FPPGOpenFree}, its abelianization $N_\ab$ is a free $\BZ_p$-module of the same rank, and \eqref{FreeSubRep1} follows for~$a$.

In the case $a\not\in N$ consider the intermediate subgroup $H:=N\langle a\rangle$. By Proposition \ref{FPPGOpenFree} this is again a finitely generated free pro-$p$-group, and the quotient $H/N$ is cyclic of order $p^r$ for some integer $r\ge1$. Applying Proposition \ref{FPPGCyclicQuot} to $N\triangleleft H$, choose generators $x_1,\ldots,x_n$ of~$H$ such that $x_2,\ldots,x_n\in N$. Then each of the residue classes $[a]$ and $[x_1]$ generates the quotient $G/N$, so after replacing $x_1$ by a suitable power we may assume that $[a]=[x_1]$.
Now $N$ is a free pro-$p$-group on the generators $x_1^{p^r}$ and $x_1^ix_jx_1^{-i}$ for all $0\le i<p^r$ and $2\le j\le n$. The residue classes of these elements thus form a basis of $N_\ab$ over~$\BZ_p$. Moreover, conjugation by $x_1$ fixes the residue class of $x_1^{p^r}$ and maps $x_1^ix_jx_1^{-i}$ to $x_1^{i+1}x_jx_1^{-i-1}$. In the case $i=p^r-1$ this image is equal to $x_1^{p^r}x_jx_1^{-p^r}$, whose residue class in $N_\ab$ coincides with that of~$x_j$, because $x_1^{p^r}$ lies in~$N$. Thus conjugation by $x_1$ transitively permutes the set of basis elements $\{[x_1^ix_jx_1^{-i}] \mid 0\le i<p^r\}$ for each $2\le j\le n$. Since $p^r>1$, this implies that the trace of this operator is equal to~$1$, proving \eqref{FreeSubRep1} for $a$. This finishes the proof of Theorem \ref{FreeSubRep}.
\end{Proof}

\medskip
Now let $G$ be a finitely generated $\sigma$-pro-$p$-group whose underlying pro-$p$-group is free
and set
\UseTheoremCounterForNextEquation
\begin{equation}\label{dGepsilon}
\begin{array}{l}
d^+_G\ :=\ \dim_{\BF_p}(G/\Fratt(G))^+, \\[3pt]
d^-_G\ :=\ \dim_{\BF_p}(G/\Fratt(G))^-.
\end{array}
\end{equation}

\begin{Prop}\label{FPPGCyclicQuotS}
\begin{enumerate}
\item[(a)] Then $G$ is a free pro-$p$-group on $d_G^+$ even and $d_G^-$ odd generators.
\item[(b)] If the subgroup $N$ in Proposition \ref{FPPGCyclicQuot} is $\sigma$-invariant, the generators $x_1,\ldots,x_n$ can be chosen such that each is even or odd.
\end{enumerate}
\end{Prop}

\begin{Proof}
Choose a basis $\bar x_1,\ldots,\bar x_n$ of $G/\Fratt(G)$ such that each $\bar x_i$ is even or odd. By Proposition \ref{SPPGExact} we can lift each $\bar x_i$ to an element $x_i\in G$ of the same parity. By Proposition \ref{NormGens} (a) the elements $x_1,\ldots,x_n$ then again generate~$G$, and as $G$ is free on $n$ generators, by Proposition \ref{AutFGCons} this implies~(a).

For (b) observe that the image of $N$ in $G/\Fratt(G)$ is a $\sigma$-invariant subspace. Thus in the proof of Proposition \ref{FPPGCyclicQuot} we can choose $\bar x_1,\ldots,\bar x_n$ such that each $\bar x_i$ is even or odd. Lifting them to elements $x_1\in G$ and $x_2,\ldots,x_n \in N$ of the same parity by Proposition \ref{SPPGExact} thus yields~(b).
\end{Proof}

\medskip
The following invariant is a kind of \emph{index} associated to~$G$:
\UseTheoremCounterForNextEquation
\begin{equation}\label{iG}
i_G\ :=\ d_G^+-d_G^-.
\end{equation}
Proposition \ref{FPPGCyclicQuotS} (a) implies that this is equal to the trace of $\sigma$ on the free $\BZ_p$-module $G_\ab$.

\begin{Prop}\label{FPPGiGSub}
For any $\sigma$-invariant open normal subgroup $N$ of $G$ we have 
$$i_N-1\ =\ |(G:N)^+|\cdot(i_G-1).$$
\end{Prop}

\begin{Proof}
We prove this by induction on $[G:N]$, the case $[G:N]=1$ being trivial. 

If $[G:N]=p$, we choose generators $x_1,\ldots,x_n$ of $G$ as in Propositions \ref{FPPGCyclicQuot} and \ref{FPPGCyclicQuotS} (b) with $r=1$. so that the residue classes of $x_1^p$ and $x_1^ix_jx_1^{-i}$ for all $0\le i<p$ and $2\le j\le n$ form a basis of $N/\Fratt(N)$.

Suppose first that $x_1$ is even. Then $x_1^p$ is even, and for any $0\le i<p$ and $2\le j\le n$ with ${}^\sigma x_j=x_j^\epsilon$ the computation 
$${}^\sigma(x_1^ix_jx_1^{-i})\ =\ {}^\sigma x_1^i\kern2pt{}^\sigma\kern-1pt x_j\kern1pt{}^\sigma\kern-1pt x_1^{-i}\ =\ x_1^i\kern1pt x_j^\epsilon\kern1pt x_1^{-i}\ =\ (x_1^i\kern1pt x_j\kern1pt x_1^{-i})^\epsilon$$
shows that $x_1^i\kern1pt x_j\kern1pt x_1^{-i}$ has the same parity as~$x_j$. Since $d^+_G-1$ of the elements $x_2,\ldots,x_n$ are even and $d_G^-$ of them are odd, it follows that $d_N^+ = 1+p(d^+_G-1)$ and $d_N^- = pd^-_G$, and thus
$$i_N-1\ =\ d^+_N-1-d^-_N\ =\ p\cdot(d^+_G-1-d^-_G)\ =\ |(G:N)^+|\cdot(i_G-1),$$
as desired.

Suppose now that $x_1$ is odd. Then $x_1^p$ is odd; hence $d^+_G$ of the elements $x_1^p,x_2,\ldots,x_n$ are even and $d^-_G$ of them are odd. For any $0< i<p$ and $2\le j\le n$ we compute 
$${}^\sigma(x_1^ix_jx_1^{-i})
\ =\ {}^\sigma x_1^i\kern2pt{}^\sigma\kern-1pt x_j\kern1pt{}^\sigma\kern-1pt x_1^{-i}
\ =\ x_1^{-i}\kern1pt\kern-1pt x_j^{\pm1}\kern1pt x_1^i
\ =\ x_1^{-p}(x_1^{p-i}\kern1pt\kern-1pt x_j\kern1pt x_1^{i-p})^{\pm1} x_1^p.$$
Since $x_1^p\in N$, this shows that ${}^\sigma(x_1^ix_jx_1^{-i})$ is congruent to $(x_1^{p-i}\kern1pt\kern-1pt x_j\kern1pt x_1^{i-p})^{\pm1}$ modulo $\Fratt(N)$. Since $i\not\equiv p-i$ modulo $(p)$, it follows that the residue classes of $x_1^ix_jx_1^{-i}$ for all $0< i<p$ and $2\le j\le n$ are pairwise interchanged up to sign by~$\sigma$. For any such pair $z,{}^\sigma\kern-1pt  z$ the element $z+{}^\sigma\kern-1pt  z$ is even and the element $z-{}^\sigma\kern-1pt  z$ is odd. Thus the $\BF_p$-subspace of $N/\Fratt(N)$ that is generated by the residue classes $x_1^ix_jx_1^{-i}$ for all $0< i<p$ and $2\le j\le n$ has even and odd part of the same dimension 
$\frac{(p-1)(n-1)}{2}$. Together this shows that $d_N^\epsilon = d_G^\epsilon + \frac{(p-1)(n-1)}{2}$ for every~$\epsilon$, which in turn implies that $i_N=d_N^+-d_N^- = d_G^+-d_G^- = i_G$. As $|(G:N)^+|=1$ in this case, the desired formula again follows.

Finally, if $[G:N]>p$, using Proposition \ref{SPGIndLem} we choose a $\sigma$-invariant normal subgroup $M\triangleleft G$ of index~$p$ containing~$N$. Then the proposition is already proved for $M\triangleleft G$, and by the induction hypothesis it holds for $N\triangleleft M$. Since $|(G/N)^+| = |(G/M)^+|\cdot|(M/N)^+|$ by Proposition \ref{SPGExact}, the proposition also follows for $N\triangleleft G$, and we are done.
\end{Proof}

\medskip
Now observe that the conjugation of $G$ on itself and the action of $\sigma$ induce a left action of the semi-direct product $G\rtimes\{1,\sigma\}$ on~$G$. For any $\sigma$-invariant open normal subgroup $N\triangleleft G$ this action induces an action on $N$, and the resulting action on $N_\ab$ factors through the group $\Delta := (G/N)\rtimes\{1,\sigma\}$.
On the other hand, for any $\epsilon\in\{\pm1\}$ consider the left $\BQ_p[\Delta]$-submodule $V_\epsilon := \BQ_p[\Delta](\sigma+\epsilon)$ of the group ring $\BQ_p[\Delta]$ and abbreviate $V_+:=V_{+1}$ and $V_-:=V_{-1}$. Also let $\BQ_p$ carry the trivial representation of~$\Delta$.

\begin{Thm}\label{FreeSubRepSigma}
There exists an isomorphism of $\BQ_p[\Delta]$-modules
$$(N_\ab\otimes_{\BZ_p}\BQ_p) \oplus V_+\ \cong\ \BQ_p \oplus V_+^{\oplus d_G^+} \oplus V_-^{\oplus d_G^-}.$$
\end{Thm}


\begin{Proof}
Let $\chi_N$, $\chi_0$, $\chi_\epsilon$ denote the respective character of $\Delta$ on $N_\ab\otimes_{\BZ_p}\BQ_p$, $\BQ_p$, $V_\epsilon$. Then it suffices to show that 
$\chi_N = \chi_0+(d_G^+-1)\chi_++d_G^-\chi_-$. By Propositions \ref{SPGSemidirect} and \ref{SPPGExact} every element of $\Delta\setminus G/N$ is conjugate to an element of the form $[a]\sigma$ for some $a\in G^+$. Thus the desired equality follows from the character values at $\delta\in\Delta$ in the following table:
\UseTheoremCounterForNextEquation
\begin{equation}\label{FreeSubRepSigma1}
\begin{array}{|c|c||c|c|c|}
\hline
{\Large\mathstrut}  \delta & \hbox{conditions} & \chi_N(\delta) & \chi_0(\delta) & \chi_\epsilon(\delta) \\
\hline\hline
{\Large\mathstrut} 1 & \hbox{none} & 1 + [G:N]\cdot(d_G^+-1+d_G^-) & 1 & [G:N] \\
\hline
{\Large\mathstrut} [a] & a\in G\setminus N & 1 & 1 & 0 \\
\hline
{\large\mathstrut} \sigma & \hbox{none} &  1+ |(G/N)^+|\cdot(d_G^+-1-d_G^-) & 1 & \epsilon\cdot|(G/N)^+| \\
\hline
{\Large\mathstrut} [a]\sigma & a\in G^+\setminus N & 1 & 1 & 0 \\
\hline\end{array}
\end{equation}

Here the values $\chi_0(\delta)=1$ are obvious. Next, each $V_\epsilon$ is a free module of rank $1$ over $\BQ_p[G/N]$ with basis $\sigma+\epsilon$. This yields the upper two values of $\chi_\epsilon(\delta)$. Also the elements $[x](\sigma+\epsilon)$ for all residue classes $[x]\in G/N$ form a basis of $V_\epsilon$ over~$\BQ_p$. The computation 
$$\sigma[x](\sigma+\epsilon)
\ =\ [{}^\sigma x] \sigma(\sigma+\epsilon)
\ =\ \epsilon [{}^\sigma x] (\sigma+\epsilon)$$ 
thus shows that $\chi_\epsilon(\sigma) = \epsilon\cdot\trace(\sigma|\BQ_p[G/N])$. Since $\sigma$ acts on the basis $G/N$ of $\BQ_p[G/N]$ by permutation, its trace is simply the number of fixed points, yielding the value $\chi_\epsilon(\sigma) = \epsilon\cdot|(G/N)^+|$. Moreover, for any $a\in G^+\setminus N$ the computation 
$$[a]\sigma[x](\sigma+\epsilon)
\ =\ [a] \epsilon [{}^\sigma x] (\sigma+\epsilon)
\ =\ \epsilon [a\kern2pt{}^\sigma x] (\sigma+\epsilon)$$
shows that $[a]\sigma$ maps each basis vector to $\pm$ a different basis vector. This implies the value $\chi_\epsilon([a]\sigma) = 0$, finishing the computation of~$\chi_\epsilon$.

Since $d_G=d_G^++d_G^-$, the upper two values of $\chi_N(\delta)$ are already known from \eqref{FreeSubRep1}. Next, Proposition \ref{FPPGCyclicQuotS} (a) implies that $N_\ab$ is a free $\BZ_p$-module on $d_N^+$ even and $d_N^-$ odd elements. This implies that $\chi_N(\sigma) = d_N^+-d_N^- = i_N$, and so the desired value for $\chi_N(\sigma)$ results from Proposition \ref{FPPGiGSub} and the equation $i_G=d_G^+-d_G^-$.

Finally take any $a\in G^+\setminus N$ and consider the intermediate subgroup $H:=N\langle a\rangle$. By Proposition \ref{FPPGOpenFree} this is again a finitely generated free pro-$p$-group, and the quotient $H/N$ is cyclic of order $p^r$ for some integer $r\ge1$. Using Propositions \ref{FPPGCyclicQuot} and \ref{FPPGCyclicQuotS} (b), choose free generators $x_1,\ldots,x_n$ of $H$ such that each $x_j$ is even or odd and that $N$ is free on the generators $x_1^{p^r}$ and $x_1^ix_jx_1^{-i}$ for all $0\le i<p^r$ and $2\le j\le n$. The residue classes of these elements then form a basis of $N_\ab$ as a free $\BZ_p$-module.
Now each of the residue classes $[a]$ and $[x_1]$ generates the quotient $H/N$, so after replacing $x_1$ by a suitable power we may assume that $[a]=[x_1]$. Also, since $a$ is even, the element $x_1$ must be even. The action of $[a]\sigma=[x_1]\sigma$ therefore fixes the basis element $[x_1^{p^r}]$ of~$N_\ab$. On the other hand, for any $0\le i<p^r$ and $2\le j\le n$ with ${}^\sigma x_j=x_j^\epsilon$ we compute 
$${}^{x_1\sigma}(x_1^ix_jx_1^{-i})\ =\ x_1\kern2pt{}^\sigma x_1^i\kern2pt{}^\sigma\kern-1pt x_j\kern1pt{}^\sigma\kern-1pt x_1^{-i}x_1^{-1}\ =\ x_1^{i+1}\kern1pt x_j^\epsilon\kern1pt x_1^{-i-1}\ =\ (x_1^{i+1}\kern1pt x_j\kern1pt x_1^{-i-1})^\epsilon.$$
In the case $i=p^r-1$ this is equal to $(x_1^{p^r}x_jx_1^{-p^r})^\epsilon$, whose residue class in $N_\ab$ coincides with that of~$x_j^\epsilon$, because $x_1^{p^r}$ lies in~$N$. Thus the action of $[a]\sigma=[x_1]\sigma$ transitively permutes each set of basis elements $\{[x_1^ix_jx_1^{-i}] \mid 0\le i<p^r\}$ up to sign. Since $p^r>1$, this implies that the trace of the total operator is equal to~$1$, yielding the desired value $\chi_N([a]\sigma)=1$. 

Having justified all entries of \eqref{FreeSubRepSigma1}, this finishes the proof of Theorem \ref{FreeSubRepSigma}.
\end{Proof}

\begin{Cor}\label{NabOpen}
The $\BQ_p[G/N]$-module $N_\ab\otimes_{\BZ_p}\BQ_p$ is generated by $d_G^+$ even and $d_G^-$ odd elements.
\end{Cor}

\begin{Proof}
The non-zero $\Delta$-invariant element $\sum_{\delta\in\Delta}\delta \in V_+$ induces an injective homomorphism of $\BQ_p[\Delta]$-modules $i\colon\BQ_p\into V_+$. The isomorphism from Theorem \ref{FreeSubRepSigma} therefore also yields an isomorphism
$$(N_\ab\otimes_{\BZ_p}\BQ_p) \oplus V_+/i(\BQ_p)\ \cong\ V_+^{\oplus d_G^+} \oplus V_-^{\oplus d_G^-}.$$
Thus there exists a surjective homomorphism of $\BQ_p[\Delta]$-modules $
V_+^{\oplus d_G^+} \oplus V_-^{\oplus d_G^-} \onto N_\ab\otimes_{\BZ_p}\BQ_p$.
Since each $V_\epsilon$ is generated as a $\BQ_p[G/N]$-module by the element $\sigma+\epsilon$ of sign~$\epsilon$, the desired assertion follows.
\end{Proof}

\section{Automorphisms of free $\sigma$-pro-$p$-groups}
\label{AutFSG}

{}From now on, for any integer $n\ge0$ we let $F_n$ denote the free pro-$p$-group with $n$ generators $x_1,\ldots,x_n$ and turn it into a $\sigma$-pro-$p$-group by requiring $x_1,\ldots,x_n$ to be \emph{odd}. Then $F_n/\Fratt(F_n)$ is an $\BF_p$-vector space of dimension~$n$ and is totally odd. For any $n\ge m\ge0$ we identify $F_m$ with a subgroup of $F_n$ 
by Fried-Jarden \cite[Lemma 20.4.9]{FriedJarden2023}.

\medskip
In this section we establish some basic properties of the group of $\sigma$-automorphisms $\Aut_\sigma(F_n)$ of~$F_n$. For any $\alpha\in\Aut_\sigma(F_n)$ we have $\alpha(x_1),\ldots,\alpha(x_n)\in F_n^-$ and their residue classes form an $\BF_p$-basis of $F_n/\Fratt(F_n)$. Conversely we have:

\begin{Prop}\label{AutFSGCons}
For any $y_1,\ldots,y_n\in F_n^-$ whose residue classes form an $\BF_p$-basis of $F_n/\Fratt(F_n)$, there exists a unique $\alpha\in\Aut_\sigma(F_n)$ with $\alpha(x_i)=y_i$ for all~$i$.
\end{Prop}

\begin{Proof}
As all $y_i$ lie in~$F_n^-$, the automorphism $\alpha$ from Proposition \ref{AutFGCons} is $\sigma$-equivariant.
\end{Proof}

\begin{Prop}\label{AutGLn}
Every automorphism $\bar\alpha$ of the $\BF_p$-vector space $F_n/\Fratt(F_n)$ can be lifted to some $\alpha\in\Aut_\sigma(F_n)$.
\end{Prop}

\begin{Proof}
Since $F_n^-\onto (F_n/\Fratt(F_n))^- = F_n/\Fratt(F_n)$ is surjective by Proposition \ref{SPPGExact}, for any $1\le i\le n$ we can choose an element $y_i\in F_n^-$ with residue class $\bar\alpha([x_i])$. Then the $\sigma$-automorphism furnished by Proposition \ref{AutFSGCons} has the desired property.
\end{Proof}

\begin{Prop}\label{AutLift}
For every $\sigma$-pro-$p$-group $G$ and any surjective $\sigma$-homomorphisms $\phi,\psi\colon \allowbreak F_n\onto G$ there exists a $\sigma$-automorphism $\alpha$ of $F_n$ such that $\phi\circ\alpha=\psi$.
\end{Prop}

\begin{Proof}
Abbreviate $N:=\ker(\phi)$. Then $\phi$ induces a surjective $\sigma$-homomorphism
\UseTheoremCounterForNextEquation
\begin{equation}\label{AutLift1}
F_n/\Fratt(F_n) \longonto F_n/\Fratt(F_n)N \;\cong\; G/\Fratt(G).
\end{equation}
In particular we have $m := \dim_{\BF_p}(G/\Fratt(G))\le \dim_{\BF_p}(F_n/\Fratt(F_n)) = n$. Thus there exists an $\BF_p$-basis $\bar y_1,\ldots,\bar y_n$ of $F_n/\Fratt(F_n)$ such that $\bar y_1,\ldots,\bar y_m$ map to a basis of $G/\Fratt(G)$ and $\bar y_{m+1},\ldots,\bar y_n$ map to zero. Since $F_n^-\onto (F_n/\Fratt(F_n))^- = F_n/\Fratt(F_n)$ is surjective by Proposition \ref{SPPGExact}, for every $1\le i\le m$ we can choose an element $y_i\in F_n^-$ with residue class~$\bar y_i$. Likewise, since $N$ surjects to the kernel $\Fratt(F_n)N/\Fratt(F_n)$ of \eqref{AutLift1}, the induced map $N^- \to (\Fratt(F_n)N/\Fratt(F_n))^- = \Fratt(F_n)N/\Fratt(F_n)$ is surjective by Proposition \ref{SPPGExact}, so for any $m< i\le n$ we can choose an element $y_i\in N^-$ with residue class~$\bar y_i$. Then the residue classes of $y_1,\ldots,y_n\in F_n^-$ form an $\BF_p$-basis of $F_n/\Fratt(F_n)$, so by Proposition \ref{AutFSGCons} there exists a $\sigma$-automorphism $\beta$ of $F_n$ such that $\beta(x_i)=y_i$ for all~$i$. By construction $\phi\circ\beta$ induces an isomorphism $F_m/\Fratt(F_m) \cong G/\Fratt(G)$ and we have $\phi(\beta(x_i))=1$ for all $m<i\le n$. In particular $\phi\circ\beta|F_m$ is surjective by Proposition \ref{NormGens}.

In the same way we can find a $\sigma$-automorphism $\gamma$ of $F_n$ such that $\psi\circ\gamma$ induces an isomorphism $F_m/\Fratt(F_m) \cong G/\Fratt(G)$ and $\psi(\gamma(x_i))=1$ for all $m<i\le n$.

Since $\phi\circ\beta|F_m$ is surjective, it also induces a surjection $F_m^-\onto G^-$ by Proposition \ref{SPPGExact}. Thus for every $1\le i\le m$ we can choose an element $z_i\in F_n^-$ with $\phi(\beta(z_i)) = \psi(\gamma(x_i))$. The residue classes of $z_1,\ldots,z_m$ then form a basis of $F_m/\Fratt(F_m)$. Setting $z_i := x_i$ for all $m<i\le n$, the residue classes of $z_1,\ldots,z_n$ thus form a basis of $F_n/\Fratt(F_n)$. By Proposition \ref{AutFSGCons} there thus exists a $\sigma$-automorphism $\delta$ of $F_n$ with $\delta(x_i)=z_i$ for all~$i$. Its construction implies that $\phi\circ\beta\circ\delta = \psi\circ\gamma$; hence the proposition holds with $\alpha = \beta\circ\delta\circ\gamma^{-1}$.
\end{Proof}

\section{Weak Schur $\sigma$-groups}
\label{WSSG}

For any tuple $\rtuple=(r_1,\ldots,r_n) \in (F_n^-)^n$ we let $N_\rtuple$ denote the normal subgroup of $F_n$ that is generated by all conjugates of $r_1,\ldots,r_n$. By construction this subgroup is $\sigma$-invariant; hence the factor group $G_\rtuple := F_n/N_\rtuple$ is a $\sigma$-pro-$p$-group defined by $n$ odd generators and $n$ odd relations $r_1,\ldots,r_n$.

\medskip
For another interpretation of this observe that by the universal property of~$F_n$, giving a tuple $\rtuple=(r_1,\ldots,r_n) \in (F_n)^n$ is equivalent to giving a homomorphism $\rho\colon F_n\to F_n$ by the relation $\rho(x_i)=r_i$. Moreover, since the generators $x_1,\ldots,x_n$ of $F_n$ are all odd, the entries of $\rtuple$ are odd if and only if $\rho$ is $\sigma$-equivariant. Giving a tuple $\rtuple \in (F_n^-)^n$ is therefore equivalent to giving a $\sigma$-homomorphism $\rho\colon F_n\to F_n$. Moreover, the subgroup $N_\rtuple$ is then simply the normal closure of the image $\rho(F_n)$.

The set of $\sigma$-homomorphisms $F_n\to F_n$ carries natural commuting left and right actions of $\alpha,\beta\in\Aut_\sigma(F_n)$ by $\rho\mapsto \alpha\circ\rho\circ\beta$. We denote the corresponding actions on the set $(F_n^-)^n$ by $\rtuple\mapsto \alpha\rtuple\beta$. Since $(\alpha\circ\rho\circ\beta)(F_n) = \alpha(\rho(F_n))$, it follows that $N_{\alpha\rtuple\beta} = \alpha(N_\rtuple)$ and that $\alpha$ induces a $\sigma$-isomorphism 
\UseTheoremCounterForNextEquation
\begin{equation}\label{PhiPsiGrIsom}
G_\rtuple \isoto G_{\alpha\rtuple\beta}.
\end{equation}
To understand these actions it helps to notice that left multiplication by $\alpha$ modifies the generators of $G$ and right multiplication by $\beta$ transforms the relations of~$G$.

\medskip
Next we identify $F_n/\Fratt(F_n)$ with the space of column vectors $\BF_p^n$ such that the images $\bar x_1,\ldots,\bar x_n$ of $x_1,\ldots,x_n$ correspond to the standard basis.
For any tuple  $\bar\rtuple=(\bar r_1,\ldots,\bar r_n) \in (\BF_p^n)^n$ we consider the subset 
\UseTheoremCounterForNextEquation
\begin{equation}\label{XrDef}
X_{\bar\rtuple}\ :=\ \bigl\{(r_1,\ldots,r_n) \in (F_n^-)^n \bigm| \forall i\colon r_i\mapsto \bar r_i \bigr\}.
\end{equation}
For any $\rtuple\in X_{\bar\rtuple}$ the image of $N_\rtuple$ in $F_n/\Fratt(F_n)\cong\BF_p^n$ is the subspace generated by $\bar r_1,\ldots,\bar r_n$ and therefore only depends on~$\bar\rtuple$. Viewing $\bar\rtuple=(\bar r_1,\ldots,\bar r_n)$ as an $n\times n$-matrix over~$\BF_p$,
the dimension of this subspace is then simply the rank of~$\bar\rtuple$. Taking the factor space it follows that the minimal number of generators of $G_\rtuple$ is
\UseTheoremCounterForNextEquation
\begin{equation}\label{dGrBarr}
d_{G_\rtuple}\ =\ \dim_{\BF_p}(G_\rtuple/\Fratt(G_\rtuple))\ =\ n - \rank(\bar\rtuple).
\end{equation}
Moreover, for any $\sigma$-automorphisms $\alpha,\beta$ of $F_n$ consider the induced automorphisms $\bar\alpha,\bar\beta \in \Aut_{\BF_p}(F_n/\Fratt(F_n)) $. Viewing $\bar\alpha, \bar\beta$ as matrices in $\GL_n(\BF_p)$, the definition of the left and right action on $(F_n^-)^n$ implies that 
\UseTheoremCounterForNextEquation
\begin{equation}\label{XrAut}
X_{\bar\alpha\kern1pt\cdot\kern1pt\bar\rtuple\kern1pt\cdot\kern1pt\bar\beta}\ =\ \alpha X_{\bar\rtuple}\, \beta.
\end{equation}
Since any two $n\times n$-matrices of the same rank can be transformed into each other through left and right multiplication by $\GL_n(\BF_p)$, and the reduction map $\Aut_\sigma(F_n)\to\GL_n(\BF_p)$ is surjective by Proposition \ref{AutGLn}, we can utilize this to reduce ourselves to the special case that $\bar\rtuple=(0,\ldots,0,\bar x_{m+1},\ldots,\bar x_n)$ for some $0\le m\le n$.

\medskip
In that case observe that for any tuple $\rtuple\in X_{\bar\rtuple}$ we have $r_i\in \Fratt(F_n)$ for $i\le m$ and $r_i\in x_i\Fratt(F_n)$ for $i>m$. Thus the residue classes of $x_1,\ldots,x_m,r_{m+1},\ldots,r_n$ form an $\BF_p$-basis of $F_n/\Fratt(F_n)$, and so $F_n$ is a free pro-$p$-group on these elements by Proposition \ref{AutFGCons}. By the universal property of $F_n$ there therefore exists a unique homomorphism ${\FEpi_\rtuple\colon F_n\to F_m}$ such that $\FEpi_\rtuple(x_i)=x_i$ for $i\le m$ and $\FEpi_\rtuple(r_i)=1$ for $i>m$. By construction this homomorphism is surjective, and since the elements $x_1,\ldots,x_m,r_{m+1},\ldots,r_n$ are all odd, it is also $\sigma$-equivariant. Moreover, for any $i\le m$ the fact that $r_i\in\Fratt(F_n)$ implies that $\FEpi_\rtuple(r_i) \in\Fratt(F_m)$. We therefore have a well-defined map
\UseTheoremCounterForNextEquation
\begin{equation}\label{ThetaDef}
\XrRedMap\colon\ X_{\bar\rtuple} \longto (\Fratt(F_m)^-)^m,\ \ 
r \longmapsto (\FEpi_\rtuple(r_1),\ldots,\FEpi_\rtuple(r_m)).
\end{equation}

\begin{Prop}\label{GrGthetar}
In the above situation
$\FEpi_\rtuple$ induces a $\sigma$-isomorphism $G_\rtuple\isoto G_{\XrRedMap(\rtuple)}$.
\end{Prop}

\begin{Proof}
Since $F_n$ is free on the generators $x_1,\ldots,x_m,r_{m+1},\ldots,r_n$ and $\FEpi_\rtuple$ sends the first $m$ of them to the free generators of $F_m$ and the last $n-m$ of them to~$1$, the homomorphism is surjective and its kernel is the normal subgroup generated by all conjugates of $r_{m+1},\ldots,r_n$. Since $N_\rtuple$ is the normal subgroup of $F_n$ that is generated by all conjugates of $r_1,\ldots,r_n$, its image is therefore the normal subgroup of $F_m$ that is generated by all conjugates of $\FEpi_\rtuple(r_1),\ldots,\FEpi_\rtuple(r_m)$. By the definition of $\XrRedMap(\rtuple)$ this means that $\FEpi_\rtuple(N_\rtuple) = N_{\XrRedMap(\rtuple)}$. As $\FEpi_\rtuple$ is surjective, it thus induces a $\sigma$-isomorphism $G_\rtuple = F_n/N_\rtuple \isoto F_m/\FEpi_\rtuple(N_\rtuple) = F_m/N_{\XrRedMap(\rtuple)} = \nobreak G_{\XrRedMap(\rtuple)}$, as desired.
\end{Proof}

\medskip
After these preliminaries we can now introduce the main object of our interest:

\begin{Prop}\label{WSSGEquiv}
For any $\sigma$-pro-$p$-group $G$ the following conditions are equivalent:
\begin{enumerate}
\item[(a)] We have $G\cong G_\rtuple$ for some $n\ge0$ and $\rtuple\in (F_n^-)^n$.
\item[(b)] We have $\dim_{\BF_p}\!H^2(G,\BF_p) \le \dim_{\BF_p}\!H^1(G,\BF_p) < \infty$ and $\sigma$ acts by $-1$ on both spaces.
\end{enumerate}
Moreover, if these conditions hold, then (a) holds with $n=d_G$.
\end{Prop}

\begin{Proof}
First suppose that $G\cong G_\rtuple$ for $\rtuple\in (F_n^-)^n$. Then by the above remarks and the $\sigma$-isomorphism \eqref{PhiPsiGrIsom}, we may without loss of generality assume that $\rtuple\in X_{\bar\rtuple}$ for $\bar\rtuple=(0,\ldots,0,\bar x_{m+1},\ldots,\bar x_n)$. Using Proposition \ref{GrGthetar} we can next replace $n$~and~$\rtuple$ by $m$~and~$\XrRedMap(\rtuple)$, after which we have 
$N_\rtuple\subset \Fratt(F_n)$. In that case the space $H^1(G,\BF_p) \cong \Hom(G/\Fratt(G),\BF_p)$ has dimension $d_G=n$ by equation \eqref{dGrBarr}, and since the generators $x_1,\ldots,x_n$ are all odd, the involution $\sigma$ acts by $-1$ on it. 
Similarly, the definition of $N_\rtuple$ implies that the residue classes of $r_1,\ldots,r_n$ generate the $\BF_p$-vector space $N/\Fratt(N)[F_n,N]$.
By the natural isomorphism \eqref{H2Cong} the space $H^2(G,\BF_p)$ therefore has dimension $\le n$, and since the relations $r_1,\ldots,r_n$ are all odd, it follows that $\sigma$ acts by $-1$ on it. Thus $G$ has the properties in (b), proving the implication (a)$\Rightarrow$(b) as well as the last statement.

\medskip
The reverse implication was proved in a similar form by Koch and Venkov in \cite[Prop.\,1, Thm.\,1]{KochVenkov1974} or \cite[Lemma\,1, Satz\,1]{KochVenkov1975}. Suppose that $G$ satisfies the conditions in~(b). Then the isomorphism \eqref{H1Cong} implies that $G/\Fratt(G)$ is an $\BF_p$-vector space of finite dimension $n:=\dim_{\BF_p}\!H^1(G,\BF_p)$ and equal to $(G/\Fratt(G))^-$. Since $G^-\to(G/\Fratt(G))^-$ is surjective by Proposition \ref{SPPGExact}, we can choose elements $a_1,\ldots,a_n\in G^-$ whose residue classes form a basis of $G/\Fratt(G)$. These elements then generate $G$ by Proposition \ref{NormGens}~(a). Next, by the universal property of $F_n$ there exists a surjective $\sigma$-homomorphism $\phi\colon F_n\onto G$ with $\phi(x_i)=a_i$ for all~$i$. Let $N$ denote its kernel. Then the conditions in (b) and the isomorphism \eqref{H2Cong} imply that $N/\Fratt(N)[G,N]$ is an $\BF_p$-vector space of dimension $\le n$ and equal to $(N/\Fratt(N)[G,N])^-$. It can therefore be generated by $n$ elements, and since $N^-\to(N/\Fratt(N)[G,N])^-$ is surjective by Proposition \ref{SPPGExact}, we can choose elements $r_1,\ldots,r_n\in N^-$ whose residue classes form a basis of $N/\Fratt(N)[G,N]$. Setting $\rtuple := (r_1,\ldots,r_n) \in (F_n^-)^n$, the subgroup $N_\rtuple$ then satisfies $N_\rtuple\Fratt(N)[G,N]=N$. By Proposition \ref{NormGens} (b) we conclude that $N_\rtuple=N$, and so $G_\rtuple = F_n/N_\rtuple\cong F_n/N\cong G$, as desired.
\end{Proof}


\begin{Def}\label{SchCons}
We call $\sigma$-pro-$p$-group with the properties in Proposition \ref{WSSGEquiv} a \emph{weak Schur $\sigma$-group}. We let $\Sch$ (pronounce \emph{Schur}) denote the set of $\sigma$-pro-$p$-groups of the form $G_\rtuple=F_n/N_\rtuple$ for all $n\ge0$ and $\rtuple\in (F_n^-)^n$, up to $\sigma$-isomorphism. For every $n\ge0$ we set
$$\begin{array}{ll}
\Sch_{\le n} \! & :=\ \bigl\{ [G]\in\Sch \bigm| d_G\le n \bigr\}
\quad\hbox{and} \\[3pt]
\Sch_n \! & :=\ \bigl\{ [G]\in\Sch \bigm| d_G=n \bigr\}.
\end{array}$$
Thus $\Sch$ can be viewed as the set of all weak Schur $\sigma$-groups up to $\sigma$-isomorphism, avoiding any set-theoretic complications.
\end{Def}


For every $n\ge0$ consider the map
\UseTheoremCounterForNextEquation
\begin{equation}\label{PinDef}
\CAN_n\colon\ (F_n^-)^n\longto\Sch,\ \ \rtuple\mapsto[G_\rtuple] = [F_n/N_\rtuple].
\end{equation}

\begin{Prop}\label{PinImage}
This map induces surjections
$$\begin{array}{cl}
(F_n^-)^n \! & \longonto\ \Sch_{\le n}
\quad\hbox{and} \\[3pt]
(\Fratt(F_n)^-)^n & \longonto\ \Sch_n.
\end{array}$$
\end{Prop}

\begin{Proof}
The first half of the proof of Proposition \ref{WSSGEquiv} shows that $\CAN_n((F_n^-)^n) \subset \Sch_{\le n}$ and $\CAN_n((\Fratt(F_n)^-)^n) = \Sch_n$ for all~$n$. It remains to show that the first inclusion is an equality. For this consider any $[G]\in\Sch_{\le n}$. By the second equality we can choose $m\le n$ and $\rtuple=(r_1,\ldots,r_m)\in (F_m^-)^m$ such that $G\cong F_m/N_\rtuple$. Set $\rtuple' :=(r_1,\ldots,r_m,x_{m+1},\ldots,x_n) \in (F_n^-)^n$. Then $N_{\rtuple'}<F_n$ contains the normal subgroup generated by all conjugates of $x_{m+1},\ldots,x_n$, whose factor group is isomorphic to~$F_m$. The image of $N_{\rtuple'}$ in $F_m$ is then simply the group~$N_\rtuple$. Thus $G\cong F_m/N_\rtuple \cong F_n/N_{\rtuple'}$, and so $[G] = \CAN_n(\rtuple') \in \CAN_n((F_n^-)^n)$, as desired.
\end{Proof}

\section{The topology}
\label{SchTop}

For any $i\ge2$ and any weak Schur $\sigma$-group $G$ we consider the subset
\UseTheoremCounterForNextEquation
\begin{equation}\label{UnGDef}
U_{i,G}\ :=\ \bigl\{ [H]\in\Sch \bigm| H/D_i(H) \cong G/D_i(G) \bigr\}.
\end{equation}

\begin{Prop}\label{TopDef}
\begin{enumerate}
\item[(a)] These sets for all $\,i$ and $G$ are a basis for a unique topology on~$\Sch$.
\item[(b)] In this topology each set $U_{i,G}$ is open and closed.
\item[(c)] For all $n\ge0$ the subsets $\Sch_{\le n}$ and $\Sch_n$ are open and closed in~$\Sch$.
\item[(d)] The topology is Hausdorff and totally disconnected.
\end{enumerate}
\end{Prop}

\begin{Proof}
For any fixed $i\ge2$, the subsets $U_{i,G}$ are the equivalence classes for the equivalence relation on~$\Sch$ that is defined by
$$[G]\sim_i [H]\ \Longleftrightarrow\ G/D_i(G) \cong H/D_i(H).$$
In particular their union is~$\Sch$. Moreover, for any $i'\ge i$ the relation $\sim_{i'}$ is a refinement of~$\sim_i$; hence every equivalence class for $\sim_i$ is a union of equivalence classes for~$\sim_{i'}$. Thus the intersection of any two such subsets is a union of other such subsets, proving (a). 

The description of the sets $U_{i,G}$ as equivalence classes also shows that the complement of each is a union of other such subsets. Since the sets are open by construction, they are therefore also closed, proving (b). 

Next, for any $[G]\in\Sch$ the minimal number of generators $d_G = \dim_{\BF_p}(G/D_2(G))$ of $G$ depends only on the equivalence class modulo $\sim_2$. For each $n\ge0$ the subsets $\Sch_{\le n}$ and $\Sch_n$ as well as their complements are therefore unions of equivalence classes for~$\sim_2$, proving (c).

Finally, for any fixed $[G]\in\Sch$ the sets $U_{i,G}$ for all $i\ge2$ form a basis of neighborhoods of $[G]$ which are both open and closed. To prove (d) it thus suffices to show that $\bigcap_i U_{i,G} = \{[G]\}$. Take any $[H]$ in this intersection, and for every $i\ge2$ let $S_i$ denote the set of $\sigma$-isomorphisms $G/D_i(G) \isoto H/D_i(H)$. Each isomorphism in $S_{i+1}$ then induces an isomorphism in~$S_i$, yielding a map $S_{i+1}\to S_i$. As the sets $S_i$ are all non-empty by assumption, and any filtered inverse limit of non-empty finite sets is non-empty, there exists a system of $\sigma$-isomorphisms $G/D_i(G) \isoto H/D_i(H)$ that is compatible with the projections $G/D_{i+1}(G)\onto G/D_i(G)$ and $H/D_{i+1}(H)\onto H/D_i(H)$. This system yields a $\sigma$-isomorphism $G \cong \invlim_i G/D_i(G) \isoto \invlim_i H/D_i(H) \cong H$. Thus $[G]=[H]$, and we are done.
\end{Proof}

\begin{Prop}\label{TopProp}
For every $n\ge0$ we have:
\begin{enumerate}
\item[(a)] The map $\CAN_n\colon (F_n^-)^n\to\Sch$ from \eqref{PinDef} is continuous.
\item[(b)] The final topology on its image $\Sch_{\le n}$ is the subspace topology induced from~$\Sch$.
\item[(c)] The subsets $\Sch_{\le n}$ and $\Sch_n$ are compact.
\end{enumerate}
\end{Prop}

\begin{Proof}
Take any tuple $\rtuple\in \CAN_n^{-1}(U_{i,G})$. By the definition of $\CAN_n$ and $U_{i,G}$ we then have $\sigma$-isomorphisms $F_n/N_\rtuple D_i(F_n) \cong G_\rtuple/D_i(G_\rtuple) \cong G/D_i(G)$. As the subgroup $N_\rtuple D_i(F_n)$ depends only on $\rtuple$ modulo $D_i(F_n)$, it follows that the whole neighborhood $(F_n^-)^n \cap \rtuple{\cdot} D_i(F_n)^n$ of $\rtuple$ is contained in $\CAN_n^{-1}(U_{i,G})$. Thus $\CAN_n^{-1}(U_{i,G})$ is open, proving (a).

Next endow $\Sch_{\le n}$ with the subspace topology induced from~$\Sch$, so that $\CAN_n$ induces a continuous surjective map $(F_n^-)^n \onto \Sch_{\le n}$ by Proposition \ref{PinImage}. Here $(F_n^-)^n$ is a closed subset of the compact space $(F_n)^n$ and is therefore itself compact. As $\Sch_{\le n}$ is Hausdorff by Proposition \ref{TopDef} (d), it follows that the map is closed. By elementary set-theoretic topology it follows that a subset $U\subset\Sch_{\le n}$ is open if and only if its inverse image $\CAN_n^{-1}(U)$ is open, proving (b).

Finally, since $(F_n^-)^n$ is compact, so is its image $\Sch_{\le n}$, and so is therefore also the closed subset $\Sch_n\subset\Sch_{\le n}$, proving (c).
\end{Proof}


\begin{Rem}\label{TopProfinite}\rm
Proposition \ref{TopDef} shows that $\Sch$ is the infinite disjoint union of the open and closed subsets~$\Sch_n$. Combined with Proposition \ref{TopProp} we see that each
$\Sch_n$ is a totally disconnected compact Hausdorff space and therefore a profinite space.
In fact, one can realize it explicitly as an inverse limit of finite discrete sets, as follows.

For every $i\ge2$ let $\Sch_n^{(i)}$ denote the set of all $\sigma$-pro-$p$-groups of the form $G_\rtuple/D_i(G_\rtuple)$ for some $n\ge0$ and $r\in (\Fratt(F_n)^-)^n$, up to $\sigma$-isomorphism. Then we have natural maps 
$$\begin{array}{rr}
\Sch_n\to\Sch_n^{(i)}, & [G_\rtuple] \mapsto [G_\rtuple/D_i(G_\rtuple)]  \rlap{\quad\hbox{and}} \\[3pt]
\Sch_n^{(i+1)}\to\Sch_n^{(i)}, &  [G_\rtuple/D_{i+1}(G_\rtuple)] \mapsto [G_\rtuple/D_i(G_\rtuple)]\rlap{.}
\end{array}$$
As in the proofs of Propositions \ref{TopDef} and \ref{TopProp} one can show that these induce a homeomorphism 
\UseTheoremCounterForNextEquation
\begin{equation}\label{SchInvLim}
\Sch_n \isoto \invlim_i \Sch_n^{(i)}.
\end{equation}
In particular this statement includes the following fact, which can also be proved independently by direct computation:
\end{Rem}

\begin{Prop}\label{SchurInvLimGroup}
For every $i\ge2$ consider a weak Schur $\sigma$-group $G_i$ and a $\sigma$-isomorphism $G_{i+1}/D_i(G_{i+1}) \isoto G_i/D_i(G_i)$. Then $G := \invlim_i G_i/D_i(G_i)$ is a weak Schur $\sigma$-group.
\end{Prop}

\medskip
Finally every finite weak Schur $\sigma$-group defines an isolated point in~$\Sch$:

\begin{Prop}\label{FinSchurOpenClosed}
\begin{enumerate}
\item[(a)] For any finite weak Schur $\sigma$-group $G$, the subset $\{[G]\}$ is open and closed in~$\Sch$.
\item[(b)] The subset $\{[G]\in\Sch\mid G\ \hbox{is finite}\}$ is open in~$\Sch$.
\end{enumerate}
\end{Prop}

\begin{Proof}
In (a) the subset is closed, because $\Sch$ is Hausdorff. To show that it is open, by Proposition \ref{WSSGEquiv} we can choose a $\sigma$-isomorphism $G\cong F_n/N_\rtuple$ with $n=d_G$ and $\rtuple\in(F_n^-)^n$. Since $G$ is finite, the subgroup $N_\rtuple$ is then open in $F_n$ and therefore finitely generated. Thus its Frattini subgroup $\Fratt(N_\rtuple)$ is open and therefore contains $D_i(F_n)$ for some $i\ge2$. 
In particular we then have $D_i(F_n)\subset N_\rtuple$ and hence $D_i(G)=1$.

We claim that $U_{i,G} = \{[G]\}$. To prove this consider any 
weak Schur $\sigma$-group $H$ with $H/D_i(H) \cong G/D_i(G)$. Since $i\ge2$ we then have $d_H = d_G \le n$, so by Proposition \ref{WSSGEquiv} again there exists a $\sigma$-isomorphism $H\cong F_n/N_{\rtuple'}$ for some $\rtuple'\in(F_n^-)^n$. This induces a surjective $\sigma$-homomorphism
$$F_n \longonto 
F_n/N_{\rtuple'}D_i(F_n)\ \cong\ 
H/D_i(H)\ \cong\ 
G/D_i(G)\ \cong\ 
G\ \cong\ 
F_n/N_\rtuple.$$
By Proposition \ref{AutLift} the composite homomorphism lifts to a $\sigma$-automorphism $\alpha$ of~$F_n$, which therefore satisfies $\alpha(N_{\rtuple'}D_i(F_n)) = N_\rtuple$. Since $\alpha(D_i(F_n)) = D_i(F_n) \subset \Fratt(N_\rtuple)$, this implies that $\alpha(N_{\rtuple'})\Fratt(N_\rtuple) = N_\rtuple$. By Proposition \ref{NormGens} (a) we thus have $\alpha(N_{\rtuple'}) = N_\rtuple$; hence $\alpha$ induces a $\sigma$-isomorphism $H\isoto G$, proving the claim.

By the definition of the topology on $\Sch$ the claim implies (a). Finally (b) is a direct consequence of~(a).
\end{Proof}

\begin{Ex}\label{Sch1}\rm
Since $F_1\cong\BZ_p$, the set $\Sch_{\le 1}$ consists of the $\sigma$-isomorphism classes $[\BZ/p^j\BZ]$ for all $j\ge0$ together with the $\sigma$-isomorphism class~$[\BZ_p]$. Of these the isomorphism class of the trivial group lies in $\Sch_0$ and the others in $\Sch_1$. Moreover $[\BZ/p^j\BZ]$ is open and closed by Proposition \ref{FinSchurOpenClosed} (a), and the sets
$$U_{i,\BZ_p}\ =\ \{\BZ_p\} \cup \{ [\BZ/p^j\BZ]\mid j\ge i\}$$
form a neighborhood base of the point $[\BZ_p]$. Thus $\Sch_{\le1}$ is homeomorphic to the one-point compactification of the discrete set $\BZ^{\ge0}$. 
\end{Ex}

\section{The probability measure}
\label{ProbMeas}

For every $n\ge0$ let $\tilde\mu_n$ denote the product measure on $(F_n^-)^n$ obtained from the Haar measure $\mu_{F_n^-}$ on each factor. 
Since the map $\CAN_n\colon (F_n^-)^n\to\Sch$ from \eqref{PinDef} is continuous by Proposition \ref{TopProp}, we can define the pushforward measure $\mu_n := \CAN_{n*}\tilde\mu_n$. By construction this is a probability measure on the Borel $\sigma$-algebra of $\Sch$, and by Proposition \ref{PinImage} it is supported on $\Sch_{\le n}$. 
Analyzing it further requires several steps.

\begin{Prop}\label{MeasLeftRightInv}
The measure $\tilde\mu_n$ is invariant under the left and right action of $\Aut_\sigma(F_n)$ from Section \ref{WSSG}.
\end{Prop}

\begin{Proof}
The bijection between $(F_n^-)^n$ and the set of $\sigma$-homomorphisms $F_n\to F_n$ comes from a bijection between $((F_n/D_i(F_n))^-)^n$ and the set of $\sigma$-homomorphisms $F_n\to F_n/D_i(F_n)$ for every $i\ge1$. Likewise, the left and right actions of $\Aut_\sigma(F_n)$ are induced by commuting left and right actions on the latter sets. As the set $((F_n/D_i(F_n))^-)^n$ is finite, the counting measure on it is invariant under these actions. The proposition thus follows from the characterization of the Haar measure on~$F_n^-$.
\end{Proof}

\begin{Lem}\label{XrMeas}
For any tuple $\bar\rtuple \in (\BF_p^n)^n$ the set $X_{\bar\rtuple}$ from \eqref{XrDef} has measure $\tilde\mu_n(X_{\bar\rtuple}) = p^{-n^2}$.
In particular we have $\tilde\mu_n((\Fratt(F_n)^-)^n) = p^{-n^2}$.
\end{Lem}

\begin{Proof}
Since $F_n/\Fratt(F_n)$ is totally odd of cardinality $p^n$, Proposition \ref{SPGSubgroup} (a) implies that the subset $F_n^-\cap g\Fratt(F_n)$ of $F_n^-$ has Haar measure $p^{-n}$ for every $g\in F_n$. As $X_{\bar\rtuple}$ and in particular $(\Fratt(F_n)^-)^n$ is a product of subsets of this form by \eqref{XrDef}, the lemma follows.
\end{Proof}


\begin{Lem}\label{ThetaStarMeas}
Consider the tuple $\bar\rtuple=(0,\ldots,0,\bar x_{m+1},\ldots,\bar x_n)$ for some $0\le m\le n$. Then the map $\XrRedMap\colon X_{\bar\rtuple} \to (\Fratt(F_m)^-)^m$ from \eqref{ThetaDef} is continuous and satisfies
$$\XrRedMap_*(\tilde\mu_n|X_{\bar\rtuple})\ =\ p^{m^2-n^2}\cdot\tilde\mu_m|((\Fratt(F_m)^-)^m.$$
\end{Lem}

\begin{Proof}
For any $i\ge2$ abbreviate $F_{n,i} := F_n/D_i(F_n)$ and $F_{m,i} := F_m/D_i(F_m)$. Then the fact that $D_i(F_n)$ is contained in $D_2(F_n) = \Fratt(F_n)$ implies that we have a natural surjection $F_{n,i} \onto F_n/\Fratt(F_n) \cong\nobreak \BF_p^n$.
By \eqref{XrDef} the residue classes modulo $D_i(F_n)$ of the elements of $X_{\bar\rtuple}$ form the set 
$$X_{\bar\rtuple,i}\ :=\ \bigl\{\rtuple_i = (r_{1,i},\ldots,r_{n,i}) \in (F_{n,i}^-)^n \bigm| \forall i\colon r_i\mapsto \bar r_i \bigr\}.$$
By the choice of $\bar\rtuple$ this subset is contained in $(\Fratt(F_{n,i})^-)^m\times (F_{n,i}^-)^{n-m}$.

For each $\rtuple\in X_{\bar\rtuple}$, by the functoriality of the Zassenhaus filtration the surjective $\sigma$-homomorphism $\FEpi_\rtuple\colon F_n\onto F_m$ induces a surjective $\sigma$-homomorphism $\FEpi_{\rtuple,i}\colon F_{n,i}\onto F_{m,i}$. The construction of $\FEpi_\rtuple$ in Section \ref{WSSG} implies that $\FEpi_{\rtuple,i}$ depends only on the image of $\rtuple$ in $X_{\bar\rtuple,i}$. The same therefore also holds for the image in $(F_{m,i}^-)^m$ of the tuple $\XrRedMap(\rtuple) = (\FEpi_\rtuple(r_1),\ldots,\FEpi_\rtuple(r_m))$ from \eqref{ThetaDef}. 
We thus get a well-defined map
$$\XrRedMap_i\colon\ X_{\bar\rtuple,i} \longto (\Fratt(F_{m,i})^-)^m,\ \ 
[\rtuple] \longmapsto ([\FEpi_\rtuple(r_1)],\ldots,[\FEpi_\rtuple(r_m)]).$$
As the source and the target of the maps $\Phi_i$ are discrete finite sets, and $\XrRedMap$ is their inverse limit over~$i$, it follows that $\XrRedMap$ is continuous.

We claim that all fibers of $\XrRedMap_i$ have the same cardinality. As the $\sigma$-homomorphism $\FEpi_\rtuple$ does not depend on the first $m$ entries of~$\rtuple$, it suffices to show that for every surjective $\sigma$-homomorphism $F_{n,i}\onto F_{m,i}$, the fibers of the induced map $\Fratt(F_{n,i})^- \to \Fratt(F_{m,i})^-$ have the same cardinality. But that is a direct consequence of Proposition \ref{SPGExact}, proving the claim.

By the construction of the Haar measure on $F_n^-$ and $F_m^-$, the claim and the fact that $\XrRedMap$ is the inverse limit of the maps $\XrRedMap_i$ imply that $\XrRedMap_*(\tilde\mu_n|X_{\bar\rtuple}) = c \cdot\tilde\mu_m|((\Fratt(F_m)^-)^m$ for some constant $c>0$. Computing the total measure on both sides with Lemma \ref{XrMeas} we find that
$$p^{-n^2}\ =\ \tilde\mu_n(X_{\bar\rtuple})\ =\ c \cdot\tilde\mu_m(((\Fratt(F_m)^-)^m)
\ =\ c\cdot p^{-m^2}.$$
Thus $c=p^{m^2-n^2}$ and we are done.
\end{Proof}


For every integer $k\ge0$ as well as for $k=\infty$ consider the positive real number 
\UseTheoremCounterForNextEquation
\begin{equation}\label{CkDef}
C_k\ :=\ \prod_{i=1}^k (1-p^{-i}).
\end{equation}

\begin{Lem}\label{RankCount}
For any $n\ge\ell\ge k\ge0$ the number of $n\times\ell$-matrices of rank $k$ over~$\BF_p$ is
$$p^{(n+\ell-k)k}\cdot \frac{C_n C_\ell}{C_{n-k} C_{\ell-k}C_k}.$$
\end{Lem}

\begin{Proof}
See for instance \cite[Thm.\,2]{FisherAlexander1966}. 
\end{Proof}


\begin{Prop}\label{MunMum}
For any $n\ge m\ge0$ we have 
$$\mu_n|\Sch_m\ =\ \frac{C_n^2}{C_m^2 C_{n-m}} \cdot \mu_m|\Sch_m.$$
\end{Prop}

\begin{Proof}
By definition the subsets $X_{\bar\rtuple}$ for all $\bar\rtuple \in (\BF_p^n)^n$ form a finite partition of $(F_n^-)^n$ by open subsets. Thus $\mu_n := \CAN_{n*}\tilde\mu_n$ is the sum of the corresponding measures $\CAN_{n*}(\tilde\mu_n|X_{\bar\rtuple})$. By \eqref{dGrBarr} the image $\CAN_n(X_{\bar\rtuple})$ is contained in $\Sch_{n-\rank(\bar\rtuple)}$; hence $\mu_n|\Sch_m$ is the sum of the measures $\CAN_{n*}(\tilde\mu_n|X_{\bar\rtuple})|\Sch_m$ for all $\bar\rtuple$ of rank $n-m$. In particular, in the special case $n=m$ this means that $\mu_m|\Sch_m = \CAN_{m*}(\tilde\mu_m|(\Fratt(F_m)^-)^m)|\Sch_m$.

Back in the general case, the isomorphism \eqref{PhiPsiGrIsom} implies that the map $\CAN_n\colon(F_n^-)^n\to\Sch$ is invariant under the left and right action of $\Aut_\sigma(F_n)$. 
Also the measure $\tilde\mu_n$ is invariant by Proposition \ref{MeasLeftRightInv}. Since the tuples of rank $n-m$ form one orbit under this action, from \eqref{XrAut} we deduce that the measures $\CAN_{n*}(\tilde\mu_n|X_{\bar\rtuple})|\Sch_m$ coincide for all these. Using Lemma \ref{RankCount} with $(n,\ell,k) = (n,n,n-m)$ it thus follows that
$$\mu_n|\Sch_m\ =\ p^{n^2-m^2}\cdot\frac{C_n^2}{C_m^2 C_{n-m}} \cdot \CAN_{n*}(\tilde\mu_n|X_{\bar\rtuple})|\Sch_m$$
for the special tuple $\bar\rtuple=(0,\ldots,0,\bar x_{m+1},\ldots,\bar x_n)$.
For that Proposition \ref{GrGthetar} implies that $\CAN_n|X_{\bar\rtuple} = \CAN_m\circ\XrRedMap$. Using Lemma \ref{ThetaStarMeas} we therefore obtain
$$\begin{array}{rl}
\CAN_{n*}(\tilde\mu_n|X_{\bar\rtuple})|\Sch_m
&=\ 
\CAN_{m*}\XrRedMap_*(\tilde\mu_n|X_{\bar\rtuple})|\Sch_m \\[3pt]
&=\ p^{m^2-n^2}\!\cdot\CAN_{m*}(\tilde\mu_m|(\Fratt(F_m)^-)^m)|\Sch_m \\[3pt]
&=\ p^{m^2-n^2}\!\cdot\mu_m|\Sch_m.
\end{array}$$
Together this implies the desired formula.
\end{Proof}


\begin{Thm}\label{MuinfMum}
For $n\to\infty$ the measures $\mu_n$ converge to a unique probability measure $\mu_\infty$ on~$\Sch$, which for every $m\ge0$ satisfies
$$\mu_\infty|\Sch_m\ =\ \frac{C_\infty}{C_m^2} \cdot \mu_m|\Sch_m.$$
\end{Thm}

\begin{Proof}
Fix any measurable subset $A\subset\Sch$. Then for any $n\ge0$, Proposition \ref{MunMum} implies that 
$$\mu_n(A)\ =\ \sum_{m=0}^\infty \mu_n(A\cap\Sch_m)
\ =\ \sum_{m=0}^\infty \frac{C_n^2}{C_m^2 C_{n-m}} \cdot \mu_m(A\cap\Sch_m).$$
To see how this behaves for $n\to\infty$ observe first that $1\ge C_k\ge C_\infty>0$ for all $k\ge0$. Since $C_k$ converges to $C_\infty$ for $k\to\infty$, for any fixed $m$ it follows that ${C_n^2}/{C_m^2 C_{n-m}}$ converges to ${C_\infty}/{C_m^2}$ for $n\to\infty$. Also, the above inequalities imply that ${C_n^2}/{C_m^2 C_{n-m}}$ is bounded independently of $n$ and~$m$.
Furthermore, the construction of $\mu_m$ and Lemma \ref{XrMeas} imply that
$$\mu_m(A\cap\Sch_m)\ \le\ \mu_m(\Sch_m)\ =\ \tilde\mu_m((\Fratt(F_m)^-)^m)
\ =\ p^{-m^2}.$$
As the sum $\sum_{m=0}^\infty p^{-m^2}$ is absolutely convergent, together this implies that
$$\lim_{n\to\infty}\mu_n(A)
\ =\ \sum_{m=0}^\infty \frac{C_\infty}{C_m^2} \cdot \mu_m(A\cap\Sch_m).$$
By the Vitali-Hahn-Saks theorem the measures $\mu_n$ thus converge to a unique measure $\mu_\infty$ that satisfies the desired equality. Finally, in the special case $A=\Sch$ we have $\mu_n(\Sch) = \tilde\mu_n((F_n^-)^n) = 1$ by construction. In the limit we therefore obtain $\mu_\infty(\Sch)=1$; hence $\mu_\infty$ is a probability measure.
\end{Proof}

\begin{Prop}\label{MuinfOfSchn}
For every $n\ge0$ we have
$$\mu_\infty(\Sch_n)\ =\ \frac{C_\infty}{C_n^2} \cdot p^{-n^2}.$$
\end{Prop}

\begin{Proof}
For any $\rtuple\in(F_n^-)^n$ we know that $[G_\rtuple]\in\Sch_n$ if and only if $\rtuple$ lies in $(\Fratt(F_n)^-)^n$. Combining Theorem \ref{MuinfMum}, the definition of $\mu_n$, and Lemma \ref{XrMeas} we deduce that
$$\mu_\infty(\Sch_n)\ =\ \frac{C_\infty}{C_n^2} \cdot \mu_n(\Sch_n)
\ =\ \frac{C_\infty}{C_n^2} \cdot \tilde\mu_n\bigl((\Fratt(F_n)^-)^n\bigr)
\ =\ \frac{C_\infty}{C_n^2} \cdot p^{-n^2},$$
as desired.
\end{Proof}

\section{Computing probabilities}
\label{CompProb}

In this section we compute the measures of certain subsets of~$\Sch$. For this we fix an integer $n\ge0$ and a $\sigma$-invariant subgroup $D$ of~$F_n$ that is contained in $\Fratt(F_n)$ and that is invariant under $\Aut_\sigma(F_n)$. Otherwise this subgroup is arbitrary. For instance, it could be a step in the descending central or derived series or the Zassenhaus filtration or be generated by the $p^n$-th powers of all elements of~$F_n$, or be constructed by a combination of these methods.

Consider any weak Schur $\sigma$-group $G$ with $d_G=n$. Then by Proposition \ref{WSSGEquiv} there exists a surjective $\sigma$-homomorphism $\pi\colon F_n\onto G$, and by Proposition \ref{AutLift} this homomorphism is unique up to composition with an element of $\Aut_\sigma(F_n)$. Thus $\pi(D)$ is a $\sigma$-invariant open normal subgroup of $G$ that is independent of the choice of~$\pi$, and so the finite factor group $G_D := G/\pi(D)$ depends only on $G$ and~$D$. 

Let $N$ denote the kernel of $\pi\colon F_n\onto G$. Then Proposition \ref{WSSGEquiv} implies that $N$ is generated by all conjugates of $n$ odd elements of~$F_n$. The same thus holds for the subgroup $N_D := ND/D$ of $F_{n,D} := F_n/D$. In particular the minimal number of generators of $N_D$ as a normal subgroup of $F_{n,D}$ is an integer $m_{D,G}\le n$. Since $N$ is unique up to $\Aut_\sigma(F_n)$ and $D$ is invariant under $\Aut_\sigma(F_n)$, this number only depends on $D$ and~$G$. 

We are interested in the set 
\UseTheoremCounterForNextEquation
\begin{equation}\label{UDGDef}
U_{D,G}\ :=\ \bigl\{ [H]\in\Sch \bigm| H_D \cong G_D \bigr\}.
\end{equation}

If $D$ is open in~$F_n$, then it contains the subgroup $D_i(F_n)$ for some $i\ge1$. 
In that case $U_{D,G}$ is the union of the subsets $U_{i,H}$ for all $[H]\in U_{D,G}$.
As the subsets $U_{i,H}$ for all $[H]\in\Sch$ form a partition of $\Sch$ by open and closed subsets, it then follows that $U_{D,G}$ is an open and closed subset of~$\Sch$.

\begin{Prop}\label{UDGMeasureOpen}
If $D$ is open in~$F_n$, then
$$\mu_\infty(U_{D,G})\ =\ \frac{C_\infty}{C_{n-m_{D,G}}}\cdot \frac{1}{\bigl|\Aut_\sigma(G_D)\bigr|}.$$
\end{Prop}

\begin{Proof}
As a preparation let $x_{1,D},\ldots,x_{n,D} \in F_{n,D}$ denote the images of the generators $x_1,\ldots,x_n$ of~$F_n$. Since $D$ is contained in $\Fratt(F_n)$, we have a surjective $\sigma$-homomorphism $F_{n,D}\onto F_n/\Fratt(F_n) \cong \nobreak \BF_p^n$, which maps the generators $x_{1,D},\ldots,x_{n,D}$ to the basis elements $\bar x_1,\ldots,\bar x_n$.

\begin{Lem}\label{AutLiftFin}
For any $\sigma$-$p$-group $H$ and any surjective $\sigma$-homomorphisms $\phi,\psi\colon F_{n,D}\onto\nobreak H$ there exists a $\sigma$-automorphism $\alpha_D$ of $F_{n,D}$ such that $\phi\circ\alpha_D=\psi$.
\end{Lem}

\begin{Proof}
Letting $\kappa$ denote the projection $F_n\onto F_{n,D}$, by Proposition \ref{AutLift} there exists a $\sigma$-automorphism $\alpha$ of $F_n$ such that $\phi\circ\kappa\circ\alpha = \psi\circ\kappa$. Since $D$ is invariant under $\Aut_\sigma(F_n)$, this yields a $\sigma$-automorphism $\alpha_D$ of $F_{n,D}$ which satisfies $\phi\circ\alpha_D\circ\kappa = \psi\circ\kappa$. As $\kappa$ is surjective, it follows that $\phi\circ\alpha_D=\psi$.
\end{Proof}

\begin{Lem}\label{AutSFin}
We have 
$$\bigl|\Aut_\sigma(F_{n,D})\bigr|\ =\ C_n\cdot |F_{n,D}^-|^n.$$
\end{Lem}

\begin{Proof}
Any automorphism $\alpha_D\in\Aut_\sigma(F_{n,D})$ is determined by the images of the generators $\alpha_D(x_{1,D}),\ldots,\alpha_D(x_{n,D})\in F_{n,D}^-$, whose residue classes form an $\BF_p$-basis of $F_n/\Fratt(F_n)$. Conversely, consider any elements $y_{1,D},\ldots,y_{n,D}\in F_{n,D}^-$ whose residue classes form an $\BF_p$-basis of $F_n/\Fratt(F_n)$. Lifting them in any way to elements $y_1,\ldots,y_n\in F_n^-$, by Proposition \ref{AutFSGCons} there exists an automorphism $\alpha\in\Aut_\sigma(F_n)$ with $\alpha(x_j)=y_j$ for all~$j$. Since $D$ is invariant under $\Aut_\sigma(F_n)$, this yields a $\sigma$-automorphism $\alpha_D$ of $F_{n,D}$ which satisfies $\alpha_D(x_{j,D})=y_{j,D}$ for all~$j$. Thus we have a bijection between $\Aut_\sigma(F_{n,D})$ and the set of tuples $(y_{1,D},\ldots,y_{n,D})\in (F_{n,D}^-)^n$ whose residue classes form an $\BF_p$-basis of $F_n/\Fratt(F_n)$. 

As the number of $\BF_p$-bases of $F_n/\Fratt(F_n)$ is $|\!\GL_n(\BF_p)|$, using Proposition \ref{SPGExact} we deduce that
$$\bigl|\Aut_\sigma(F_{n,D})\bigr|\ =\ \bigl|\GL_n(\BF_p)\bigr|\cdot\bigl|\Fratt(F_{n,D})^-\bigr|^n.$$
On the other hand, since $F_{n,D}/\Fratt(F_{n,D}) \cong F_n/\Fratt(F_n)$ is totally odd of order~$p^n$, by Proposition \ref{SPGExact} we also have $p^n\cdot|\Fratt(F_{n,D})^-| = |F_{n,D}^-|$. Since $|\GL_n(\BF_p)| = C_n\, p^{n^2}$ by Lemma \ref{RankCount} with $(n,\ell,k) = (n,n,n)$, the desired formula follows.
\end{Proof}

\begin{Lem}\label{StabSFin}
For any $\sigma$-invariant subgroup $H$ of $\Fratt(F_{n,D})$ we have 
$$\bigl|\Stab_{\Aut_\sigma(F_{n,D})}(H)\bigr|\ =\ \bigl|\Aut_\sigma(F_{n,D}/H)\bigr| \cdot |H^-|^n.$$
\end{Lem}

\begin{Proof}
Any $\sigma$-automorphism of $F_{n,D}$ that stabilizes~$H$ induces a $\sigma$-automorphism of $F_{n,D}/H$. Conversely, consider any $\sigma$-automorphism $\bar\alpha_D$ of $F_{n,D}/H$. Then applying Lemma \ref{AutLiftFin} to the projection $F_{n,D}\onto F_{n,D}/H$ and its composite with~$\bar\alpha_D$ shows that $\bar\alpha_D$ can be lifted to a $\sigma$-automorphism $\alpha_D$ of $F_{n,D}$ that stabilizes~$H$. Thus we have a surjective homomorphism
\UseTheoremCounterForNextEquation
\begin{equation}\label{StabSFinForm}
\Stab_{\Aut_\sigma(F_{n,D})}(H) \longonto \Aut_\sigma(F_{n,D}/H).
\end{equation}
The kernel of this homomorphism consists of all $\sigma$-automorphisms $\alpha_D$ of $F_{n,D}$ that are the identity modulo~$H$. This condition means that $\alpha_D(x_{j,D}) \in F_{n,D}^-\cap x_{j,D}H$ for all $1\le j\le n$.
Conversely, for any choice of elements $y_{j,D} \in F_{n,D}^-\cap x_{j,D} H$, the fact that $H$ is contained in $\Fratt(F_{n,D})$ implies that the residue classes of $y_{1,D},\ldots,y_{n,D}$ form an $\BF_p$-basis of $F_n/\Fratt(F_n)$. As in the proof of Lemma \ref{AutSFin} there therefore exists a unique $\sigma$-automorphism $\alpha_D$ of $F_{n,D}$ such that $\alpha_D(x_{j,D}) = y_{j,D}$ for all~$j$. Thus we have a bijection between the kernel of the homomorphism \eqref{StabSFinForm} and the set of tuples $(y_{1,D},\ldots,y_{n,D})$ in $\bigtimes_{j=1}^n(F_{n,D}^-\cap x_{j,D} H)$. Since all $x_{j,D}$ are odd, Proposition \ref{SPGExact} implies that this set has cardinality $|H^-|^n$. With \eqref{StabSFinForm} the lemma follows.
\end{Proof}

\medskip
For any tuple $\rtuple_D = (r_{1,D},\ldots,r_{n,D}) \in (F_{n,D}^-)^n$ let $N_{\rtuple_D,D}$ denote the normal subgroup of $F_{n,D}$ that is generated by all conjugates of $r_{1,D},\ldots,r_{n,D}$.

\begin{Lem}\label{TuplesCount}
The number of tuples $\rtuple_D \in (F_{n,D}^-)^n$ such that $N_{\rtuple_D,D} = N_D$ is 
$$\frac{C_n}{C_{n-m_{D,G}}}\cdot |N_D^-|^n.$$
\end{Lem}

\begin{Proof}
First observe that $N_{\rtuple_D,D} \subset N_D$ if and only if all $r_{j,D}$ lie in~$N_D^-$. Granting this, by Proposition \ref{NormGens} (b) we have $N_{\rtuple_D,D} = N_D$ if and only if $N_{\rtuple_D,D}\Fratt(N_D)[F_{n,D},N_D]=\nobreak N_D$. By the definition of $N_{\rtuple_D,D}$ that is equivalent to saying that the residue classes of $r_{1,D},\ldots,r_{n,D}$ generate the $\BF_p$-vector space $N_D/\Fratt(N_D)[F_{n,D},N_D]$. Abbreviating $m := m_{D,G}$, by Proposition \ref{NormGens} (b) and the definition of $m_{D,G}$ this is an $\BF_p$-vector space of dimension~$m$. Thus the number of $n$-tuples of vectors which generate that space is the number of $m\times n$-matrices of rank $m$ over~$\BF_p$. By Lemma \ref{RankCount} with $(n,\ell,k) = (n,m,m)$ the number of such matrices is $p^{nm}C_n/C_{n-m}$. 

On the other hand, any element of $N_D/\Fratt(N_D)[F_{n,D},N_D]$ is odd and thus lifts to exactly $|(\Fratt(N_D)[F_{n,D},N_D])^-|$ elements of $N_D^-$ by Proposition \ref{SPGExact}. Since $N_D/\Fratt(N_D)[F_{n,D},N_D]$ has order $p^m$, the last formula in Proposition \ref{SPGExact} implies that this number is also equal to $p^{-m}\cdot|N_D^-|$. Thus any $n$-tuple of generators of $N_D/\Fratt(N_D)[F_{n,D},N_D]$ arises from precisely $p^{-nm}\cdot|N_D^-|^n$ tuples $(r_{1,D},\ldots,r_{n,D})\in(N_D^-)^n$.
Together the stated formula follows.
\end{Proof}

\medskip
Now we can prove Proposition \ref{UDGMeasureOpen}. Since $D$ is contained in $\Fratt(F_n)$, for every isomorphism class $[H]\in U_{D,G}$ we have $H/\Fratt(H) \cong G/\Fratt(G)$ and hence $d_H=d_G=n$. By Theorem \ref{MuinfMum} it thus suffices to compute $\mu_n(U_{D,G})$, which by the definition of $\mu_n$ is equal to $\tilde\mu_n(\CAN_n^{-1}(U_{D,G}))$. 

Next observe that $G_D\cong F_n/ND$ and that $G_{\rtuple,D}\cong F_n/N_\rtuple D$ for every $\rtuple\in(F_n^-)^n$. 
By the construction of $U_{D,G}$ and $\CAN_n$ the set $\CAN_n^{-1}(U_{D,G})$ therefore consists of all $\rtuple\in(F_n^-)^n$ such that $F_n/N_\rtuple D$ is $\sigma$-isomorphic to $F_n/ND$. As this condition depends only on the residue class of $\rtuple$ modulo~$D$, by the definition of the Haar measure on $F_n^-$ the value $\tilde\mu_n(\CAN_n^{-1}(U_{D,G}))$ is $|F_{n,D}^-|^{-n}$ times the number of residue classes of such tuples modulo~$D$. These residue classes are precisely all tuples $\rtuple_D\in(F_{n,D}^-)^n$ such that $F_{n,D}/N_{\rtuple_D,D}$ is $\sigma$-isomorphic to $F_{n,D}/N_D$.

To compute their cardinality observe that any $\sigma$-isomorphism $F_{n,D}/N_D \isoto F_{n,D}/N_{\rtuple_D,D}$ can be lifted to a $\sigma$-automorphism $\alpha_D$ of $F_{n,D}$ by Lemma \ref{AutLiftFin}. This $\sigma$-automorphism then satisfies $\alpha_D(N_D)=N_{\rtuple_D,D}$. Since everything is invariant under $\Aut_\sigma(F_{n,D})$, it follows that the cardinality in question is the number of possible conjugates $\alpha_D(N_D)$ times the number of tuples $\rtuple_D\in(F_{n,D}^-)^n$ such that $N_{\rtuple_D,D}=N_D$. On the other hand the fact that $d_G=n$ implies that $N\subset\Fratt(F_n)$ and hence $N_D\subset\Fratt(F_{n,D})$. By Lemmas \ref{AutSFin} and \ref{StabSFin} the number of conjugates $\alpha_D(N_D)$ is therefore
$$\frac{\bigl|\Aut_\sigma(F_{n,D})\bigr|}{\bigl|\Stab_{\Aut_\sigma(F_{n,D})}(N_D)\bigr|}
\ =\ \frac{C_n\cdot |F_{n,D}^-|^n}{\bigl|\Aut_\sigma(F_{n,D}/N_D)\bigr| \cdot |N_D^-|^n}.$$
Multiplied by the number of tuples from Lemma \ref{TuplesCount} the desired cardinality becomes
$$\frac{C_n^2\cdot |F_{n,D}^-|^n}{C_{n-m_{D,G}}\cdot \bigl|\Aut_\sigma(F_{n,D}/N_D)\bigr|}.$$

As explained above, by the construction of the Haar measure this shows that
$$\mu_n(U_{D,G})\ =\ \tilde\mu_n(\CAN_n^{-1}(U_{D,G}))\ =\ 
\frac{C_n^2}{C_{n-m_{D,G}}}\cdot\frac{1}{\bigl|\Aut_\sigma(F_{n,D}/N_D)\bigr|}.$$
With Theorem \ref{MuinfMum} it follows that
$$\mu_\infty(U_{D,G})\ =\ \frac{C_\infty}{C_n^2} \cdot \mu_n(U_{D,G})
\ =\ \frac{C_\infty}{C_{n-m_{D,G}}}\cdot\frac{1}{\bigl|\Aut_\sigma(F_{n,D}/N_D)\bigr|}.$$
Since $F_{n,D}/N_D \cong F_n/ND \cong G_D$, this finishes the proof of Proposition \ref{UDGMeasureOpen}.
\end{Proof}

\medskip
Next we want to apply Proposition \ref{UDGMeasureOpen} to a special case. For this we first define the \emph{Zassenhaus type} of a weak Schur $\sigma$-group.

\begin{Cons}\label{ZassTypeCons}\rm
Let $G$ be a weak Schur $\sigma$-group with $d_G=n$. Among all tuples $\rtuple\in(F_n^-)^n$ with $G\cong F_n/N_\rtuple$, select one that is \emph{minimal} in the sense of Koch \cite[Satz~1]{Koch1969}. This means that for every integer $i\ge1$, an initial segment of $\rtuple$ induces a minimal system of generators of $ND_i(F_n)/D_i(F_n)$ as a normal subgroup of $F_n/D_i(F_n)$. 
For each $1\le j\le n$ let $d_j$ be the unique integer such that $r_j\in D_{d_j}(F_n)\setminus D_{d_j+1}(F_n)$ if $r_j\not=1$, respectively $d_j:=\infty$ if $r_j=1$. Generalizing McLeman \cite[end of \S2]{McLeman2008} we call the tuple $(d_1,\ldots,d_n)$ the \emph{Zassenhaus type of~$G$.} As a consequence of \cite[Satz~1]{Koch1969} this depends only on the isomorphism class of~$G$. Also, since $d_G=n$, we necessarily have $\rtuple\in(\Fratt(F_n)^-)^n$ and hence all $d_j>1$. Moreover, by Koch and Venkov \cite[Thm.\,2]{KochVenkov1974} or \cite[Satz~2]{KochVenkov1975} each $d_j$ can only be odd or~$\infty$. The construction thus implies that $3\le d_1\le\ldots\le d_n\le\infty$.

\end{Cons}

\begin{Prop}\label{UiGMeasure}
Let $G$ be a weak Schur $\sigma$-group of Zassenhaus type $(d_1,\ldots,d_n)$. Take any integer $i\ge2$ and let $m$ be the number of $j$ such that $d_j<i$. Then
$$\mu_\infty(U_{i,G})\ =\ \frac{C_\infty}{C_{n-m}}\cdot \frac{1}{\bigl|\Aut_\sigma(G/D_i(G))\bigr|}.$$
\end{Prop}

\begin{Proof}
Setting $D=D_i(F_n)$, for any weak Schur $\sigma$-group $G$ we get $G_D \cong G/D_i(G)$ and hence $U_{D,G} = U_{i,G}$. Take a minimal tuple $\rtuple$ with $G\cong F_n/N_\rtuple$ and put $N:=N_\rtuple$. Then the minimal number $m_{D,G}$ of generators of $N_D = N_\rtuple D_i(F_n)/D_i(F_n)$ as a normal subgroup of $F_n/D_i(F_n)$ is precisely the number of $j$ such that $d_j<i$. The desired formula is thus a special case of Proposition \ref{UDGMeasureOpen}.
\end{Proof}

\medskip
We can also massage Proposition \ref{UDGMeasureOpen} into allowing weaker assumptions:

\begin{Prop}\label{UDGMeasureFin}
If $G_D$ is finite, then $U_{D,G}$ is an open and closed subset of measure
$$\mu_\infty(U_{D,G})\ =\ \frac{C_\infty}{C_{n-m_{D,G}}}\cdot \frac{1}{\bigl|\Aut_\sigma(G_D)\bigr|}.$$
\end{Prop}

\begin{Proof}
The fact that $G_D \cong F_n/ND$ is finite means that $ND$ is an open subgroup of~$F_n$. Thus $ND$ is finitely generated, and so its subgroup $\Fratt(ND)[F_n,ND]$ is again open. We can therefore pick an integer $i\ge2$ with $D_i(F_n)\subset \Fratt(ND)[F_n,ND]$. Then $D' := DD_i(F_n)$ is an open subgroup of $F_n$ that is again contained in $\Fratt(F_n)$ and is invariant under $\Aut_\sigma(F_n)$. 
Since $D_i(G)$ is contained in $ND$, we then have $ND'=NDD_i(G)=ND$. In particular this implies that $G_{D'}\cong F_n/ND'=F_n/ND\cong G_D$.

\begin{Lem}\label{UDGMeasureFinLem1}
For any normal subgroup $M<F_n$ we have $MD=ND$ if and only if $MD'=ND'$.
\end{Lem}

\begin{Proof}
The inclusion $D \subset D'$ shows the implication ``$\Rightarrow$''. To prove ``$\Leftarrow$'' assume that $MD'=ND'$. By the construction of $D'$ we then have $MDD_i(F_n) = MD' = ND' = ND$. With the assumption $D_i(F_n)\subset \Fratt(ND)[F_n,ND] \subset ND$ this implies that
$$MD \Fratt(ND)[F_n,ND]\ =\ ND.$$
By Proposition \ref{NormGens} (b) we therefore have $MD = ND$, as desired.
\end{Proof}

\begin{Lem}\label{UDGMeasureFinLem2}
For any weak Schur $\sigma$-group $H$ we have $H\kern-1pt_D\cong G\kern-1pt_D$ if and only if ${H\kern-1pt_{D'}\cong G\kern-1pt_{D'}}$.
\end{Lem}

\begin{Proof}
Since $D\subset D'\subset\Fratt(F_n)$, either of the two conditions implies that $d_H=d_G=n$. Thus we can choose a tuple $\rtuple\in(F_n^-)^n$ such that $H\cong F_n/N_\rtuple$. 
Then $H_D \cong F_n/N_\rtuple D$, and Proposition \ref{AutLift} implies that any $\sigma$-isomorphism $F_n/N_\rtuple D \cong F_n/ND$ is induced by some $\sigma$-automorphism of~$F_n$. Thus we have $H_D\cong G_D$ if and only if the subgroups $N_\rtuple D$ and $ND$ are conjugate under $\Aut_\sigma(F_n)$. The same argument with $D'$ in place of $D$ shows that $H_{D'}\cong G_{D'}$ if and only if the subgroups $N_\rtuple D'$ and $ND'$ are conjugate under $\Aut_\sigma(F_n)$.
After conjugating $N_\rtuple D$ by a suitable element of $\Aut_\sigma(F_n)$, the desired equivalence thus results from Lemma \ref{UDGMeasureFinLem1}.
\end{Proof}

\begin{Lem}\label{UDGMeasureFinLem3}
We have $m_{D,G}=m_{D',G}$.
\end{Lem}

\begin{Proof}
By definition $m_{D,G}$ is the minimal number of generators of $ND/D$ as a normal subgroup of $F_n/D$. In other words it is the smallest integer $m$ for which there exists a normal subgroup $M\triangleleft F_n$ that is generated by all conjugates of $m$ elements and which satisfies $MD/D=ND/D$. But by Lemma \ref{UDGMeasureFinLem1} the last condition is equivalent to $MD'/D'=ND'/D'$. The definition of $m_{D',G}$ thus implies the desired equality.
\end{Proof}

\medskip
Combining everything, we have $U_{D,G}=U_{D',G}$ by Lemma \ref{UDGMeasureFinLem2} and $m_{D,G}=m_{D',G}$ by Lemma \ref{UDGMeasureFinLem3}. We also have $G_D\cong G_{D'}$, and $D'$ satisfies the assumptions of Proposition~\ref{UDGMeasureOpen}. Thus the formula in Proposition \ref{UDGMeasureFin} follows.
\end{Proof}


\begin{Prop}\label{UDGMeasureFinVar}
If $G_D$ is finite and $D\subset[F_n,F_n]$, then $U_{D,G}$ is an open and closed subset of measure
$$\mu_\infty(U_{D,G})\ =\ \frac{C_\infty}{\bigl|\Aut_\sigma(G_D)\bigr|}.$$
\end{Prop}

\begin{Proof}
By assumption we have $ND\subset N[F_n,F_n]$ and hence a surjective homomorphism
$$G_D\ \cong\ F_n/ND \longonto F_n/N[F_n,F_n]\ \cong\ G_\ab.$$
Since $G_D$ is finite, it follows that $G_\ab$ is finite. As a quotient of $F_n/[F_n,F_n] \cong \BZ_p^n$ the group $G_\ab$ is therefore not defined by fewer than $n$ relations. Thus the minimal number of generators of $N[F_n,F_n]/[F_n,F_n]$ as a normal subgroup of $F_n/[F_n,F_n]$ is~$\ge n$. Since $F_n/D\onto F_n/[F_n,F_n]$ induces a surjection $ND/D\onto N[F_n,F_n]/[F_n,F_n]$, it follows that the minimal number $m_{D,G}$ of generators of $ND/D$ as a normal subgroup of $F_n/D$ is also~$\ge n$. As in any case $m_{D,G}\le n$, this shows that $m_{D,G}=n$. 

Since $C_0=1$, the desired formula is therefore a special case of Proposition \ref{UDGMeasureFin}.
\end{Proof}

\begin{Cor}\label{GAbMeasure}
For any finite abelian $p$-group $A$, the subset $\{[G]\in\Sch \mid G_\ab\cong A\}$ is open and closed of measure
$$\frac{C_\infty}{\bigl|\Aut(A)\bigr|}.$$
\end{Cor}

\begin{Proof}
Pick any weak Schur $\sigma$-group $G$ with $G_\ab\cong A$, for instance the factor group of $F_n$ defined by the relations $x_1^{p^{\nu_1}}=\ldots=x_n^{p^{\nu_n}}=1$ if $A\cong \bigoplus_{j=1}^n \BZ/p^{\nu_j}\BZ$ with all $\nu_j\ge1$. Then the group $D := [F_n,F_n]$ satisfies the assumptions of Proposition \ref{UDGMeasureFinVar} and we have $G_D\cong G_\ab \cong A$. Since $G_D$ is totally odd, this also induces an isomorphism $\Aut_\sigma(G_D)\cong\Aut(A)$. As we likewise have $H_D\cong H_\ab$ for every $[H]\in \Sch$ with $d_H=n$, the corollary follows.
\end{Proof}


\begin{Prop}\label{FinSchurMeas}
If $G$ is finite, the subset $\{[G]\} \subset \Sch$ is open and closed of measure
$$\mu_\infty(\{[G]\})\ =\ \frac{C_\infty}{\bigl|\Aut_\sigma(G)\bigr|}.$$
\end{Prop}

\begin{Proof}
This is the special case $D=1$ of Proposition \ref{UDGMeasureFinVar}.
\end{Proof}

\begin{Ex}\label{Sch1Ex}\rm 
As special cases of Proposition \ref{FinSchurMeas}, the isomorphism classes of all finite cyclic groups have measures
$$\begin{array}{rll}
\mu_\infty(\{[1]\}) &=\ C_\infty & \hbox{and} \\[3pt]
\mu_\infty(\{[\BZ/p^j\BZ]\}) &=\
\displaystyle \frac{C_\infty}{p^{j-1}(p-1)}
\ =\ \frac{C_\infty}{C_1}\cdot p^{-j} & \hbox{for all $j\ge1$.}
\end{array}$$
{}From the description of the points in $\Sch_1$ in Example \ref{Sch1} and Proposition \ref{MuinfOfSchn} we deduce that
$$\mu_\infty([\BZ_p])
\ =\ \mu_\infty(\Sch_1) - \sum_{j=1}^\infty \mu_\infty(\{[\BZ/p^j\BZ]\})
\ =\ \frac{C_\infty}{C_1^2} \cdot p^{-1} - \sum_{j=1}^\infty \frac{C_\infty}{C_1}\cdot p^{-j}$$
Since $\sum_{j=1}^\infty p^{-j} = p^{-1}/C_1$, we conclude that 
$$\mu_\infty([\BZ_p])\ =\ 0.$$
This last fact is also a consequence of Theorem \ref{StrongThm} below.
\end{Ex}

\section{Strong Schur $\sigma$-groups}
\label{SSSG}


\begin{Def}\label{StrongDef}
A weak Schur $\sigma$-group $G$ such that every open subgroup $H<G$ has finite abelianization $H_\ab$ is called a \emph{strong Schur $\sigma$-group}.
\end{Def}

Recall that according to Koch-Venkov \cite[Prop.\,1, Thm.\,1]{KochVenkov1974} or \cite[Lemma\,1, Satz\,1]{KochVenkov1975} and others, a $\sigma$-pro-$p$-group $G$ such that $\dim_{\BF_p}\!H^2(G,\BF_p) = \dim_{\BF_p}\!H^1(G,\BF_p)$, that $G_\ab$ is finite, and that $\sigma$ acts by $-1$ on $H^1(G,\BF_p)$, is called a \emph{Schur $\sigma$-group}.


\begin{Prop}\label{SchurSigmaProp}
A weak Schur $\sigma$-group $G$ such that $G_\ab$ is finite is the same as a \emph{Schur $\sigma$-group}. In particular every strong Schur $\sigma$-group is a Schur $\sigma$-group.
\end{Prop}

\begin{Proof}
First consider a weak Schur $\sigma$-group $G$ such that $G_\ab$ is finite. Then $G_\ab$ is an abelian pro-$p$-group with $n:=d_G$ generators and, being finite, it cannot be defined by fewer than $n$ relations. Thus $G$ itself cannot be defined by fewer relations, and so the inequality $\dim_{\BF_p}\!H^2(G,\BF_p) \le \dim_{\BF_p}\!H^1(G,\BF_p) = n$ in Proposition \ref{WSSGEquiv} (b) must be an equality. Therefore $G$ satisfies the above stated  conditions from \cite{KochVenkov1974,KochVenkov1975}.

Conversely, consider a $\sigma$-pro-$p$-group satisfying those conditions. Since $G_\ab$ is finite, so is $n:=d_G = \dim_{\BF_p}\!H^1(G,\BF_p)$. Moreover, since both $G$ and $G_\ab$ are defined by precisely $n$ relations, the natural homomorphism $H^2(G,\BF_p)\to H^2(G_\ab,\BF_p)$ must be an isomorphism. On the other hand, as $\sigma$ acts by $-1$ on $H^1(G,\BF_p) \cong\Hom(G_\ab,\BF_p)$, and $G_\ab$ decomposes as $G_\ab^+ \times G_\ab^-$, it follows that $G_\ab=G_\ab^-$, and so $\sigma$ acts by $-1$ on~$G_\ab$. This implies that $\sigma$ acts by $-1$ on $H^2(G_\ab,\BF_p)$; hence $G$ satisfies the conditions in \ref{WSSGEquiv} (b).
\end{Proof}

\medskip
Clearly every finite weak Schur $\sigma$-group is a strong Schur $\sigma$-group. Moreover, there exist weak Schur $\sigma$-groups which are not Schur $\sigma$-groups, for instance the group $\BZ_p$ in Example \ref{Sch1}. There also exist Schur $\sigma$-groups that are not strong:

\begin{Ex}\label{SchurNotStrong}\rm
Let $G$ be the $\sigma$-pro-$p$-group with two odd generators $x$ and~$y$ and the two odd relations $x^p=y^p=1$. Then $G_\ab$ is isomorphic to $(\BZ/p\BZ)^2$ and hence finite, so $G$ is a Schur $\sigma$-group by Proposition~\ref{SchurSigmaProp}.
Also the subgroup $H := [G,G]$ is open and we claim that its abelianization $H_\ab$ is infinite. To see this, let $\zeta$ be a primitive $p$-th root of unity and consider the matrices 
$$X := \begin{pmatrix} \zeta & 0 \\ 0 & 1 \end{pmatrix}
\quad\hbox{and}\quad
Y := \begin{pmatrix} 1 & 1 \\ 0 & \zeta \end{pmatrix}$$
over $\BQ_p(\zeta)$. As these matrices satisfy $X^p=Y^p=1$ and are contained in the pro-$p$-subgroup
$$\Bigl\{ g\in\GL_2(\BZ_p[\zeta]) \Bigm| g\equiv 
\begin{pmatrix} 1 & * \\ 0 & 1 \end{pmatrix} \bmod(\zeta-1) \Bigr\},$$
there exists a homomorphism $\phi\colon G\to\GL_2(\BZ_p[\zeta])$ with $\phi(x)=X$ and $\phi(y)=Y$. This homomorphism sends the subgroup $H$ to the group of upper triangular matrices with diagonal~$1$ and the element $[x,y]\in H$ to the matrix
$$XYX^{-1}Y^{-1} = \begin{pmatrix} 1 & 1-\zeta^{-1} \\ 0 & 1 \end{pmatrix}.$$
Thus $\phi(H)$ is abelian and infinite, and hence $H_\ab$ is infinite, as claimed.
\end{Ex}

\medskip
The main goal of this section is to prove:

\begin{Thm}\label{StrongThm}
The subset $\Sch^\strong$ of $\sigma$-isomorphisms classes of strong Schur $\sigma$-groups
is a countable intersection of open subsets of $\Sch$ and has measure~$1$.
\end{Thm}

The proof of this requires some preparation. First we fix $n\ge0$ and a $\sigma$-invariant open normal subgroup $N\triangleleft F_n$. Recall from Section~\ref{FPPG} that $N_\ab$ is a finitely generated free $\BZ_p$-module with an action of $\sigma$, on which conjugation by~$F_n$ induces a left action of the finite group $F_n/N$. Consider the associated $\BQ_p[F_n/N]$-module $V_N := N_\ab\otimes_{\BZ_p}\BQ_p$.

\begin{Lem}\label{ZarClosProp}
The tuples $(v_1,\ldots,v_n)\in(V_N^-)^n$ which do not generate $V_N$ as a $\BQ_p[F_n/N]$-module form a Zariski closed proper subset of $(V_N^-)^n$.
\end{Lem}

\begin{Proof}
The condition on the tuple means that the vectors ${}^gv_i$ for all $g\in F_n/N$ and $1\le i\le n$ do not generate the $\BQ_p$-vector space~$V_N$. Identifying $V_N$ with the space of column vectors $\BQ_p^m$ for some $m\ge0$, this is equivalent to saying that the matrix with columns ${}^gv_i$ for all $g$ and $i$ in some order has rank $<m$. That is so if and only if all $m\times m$-\allowbreak{}subdeterminants of this matrix are zero. Since $F_n/N$ acts on $V_N$ through $\BQ_p$-linear maps, these equations translate into finitely many polynomial equations in the coefficients of $v_1,\ldots,v_n$ expressed in some $\BQ_p$-basis of~$V_N^-$. Thus the set in question is Zariski closed.

On the other hand for $G:=F_n$ we have $d_G^+=0$ and $d_G^-=n$. From Corollary \ref{NabOpen} we thus know that $V_N$ can be generated as a $\BQ_p[F_n/N]$-module by $n$ elements of~$V_N^-$. Thus the set is a proper subset.
\end{Proof}

\begin{Lem}\label{DiabInfClosMeas0}
The tuples $\rtuple=(r_1,\ldots,r_n)\in(N^-)^n$ for which the group $N/[N,N]N_\rtuple$ is infinite form a closed subset of measure $0$ for the Haar measure on~$(F_n^-)^n$.
\end{Lem}

\begin{Proof}
Recall that $N_\rtuple$ is the normal subgroup of $F_n$ that is generated by all conjugates of $r_1,\ldots,r_n$. In the present case it is therefore contained in~$N$, and its image in $N_\ab$ is generated by all conjugates of the residue classes $v_1,\ldots,v_n$ of $r_1,\ldots,r_n$. In other words this image is the $\BZ_p[F_n/N]$-submodule of $N_\ab$ that is generated by $v_1,\ldots,v_n$.

On the other hand recall that $N_\ab = N/[N,N]$ is a finitely generated free $\BZ_p$-module. Thus a $\BZ_p$-submodule of $N_\ab$ has infinite index if and only if it does not generate $V_N$ over~$\BQ_p$. By Lemma \ref{ZarClosProp} the tuples in question are therefore those for which $(v_1,\ldots,v_n)$ lies in a certain Zariski closed proper subset of $(V_N^-)^n$.

By Bourbaki \cite[10.1.3.a]{BourbakiVarDiff2} this is a closed subset of measure $0$ for the Haar measure on $(V_N^-)^n$. Its intersection with $(N_\ab^-)^n$ therefore has the same properties. Since by Proposition \ref{SPGQuot}, the Haar measure on $N_\ab^-$ is just the image of the Haar measure on $N^-$, the tuples in question thus form a closed subset of measure $0$ for the Haar measure on $(N^-)^n$. Finally, since $N$ is an open subgroup of~$F_n$, the same statement holds for the Haar measure on $(F_n^-)^n$ by Proposition \ref{SPGSubgroup} (c).
\end{Proof}

\begin{Prop}\label{DiableI}
For any integer $i\ge1$, the subset
$$\bigl\{[G]\in\Sch \bigm| D_i(G)_\ab \hbox{ is infinite} \bigr\}$$
is closed of measure $0$.
\end{Prop}

\begin{Proof}
Since each $\Sch_{\le n}$ is open and closed in $\Sch$, it suffices to prove this for the intersection with $\Sch_{\le n}$. By Proposition \ref{TopProp} (b) and Theorem \ref{MuinfMum} and the construction of~$\mu_n$, it thus suffices to show that the set
$$\bigl\{r\in(F_n^-)^n \bigm| D_i(G_\rtuple)_\ab \hbox{ is infinite} \bigr\}$$
is closed of measure~$0$ for the Haar measure $\tilde\mu_n$ on $(F_n^-)^n$.
Moreover, as the subgroup $D_i(F_n)$ has finite index in~$F_n$, it suffices to show this for all tuples $\rtuple$ in a fixed residue class modulo $D_i(F_n)$. Then the $\sigma$-invariant open subgroup $N := D_i(F_n)N_\rtuple$ is independent of~$\rtuple$, and from $G_\rtuple=F_n/N_\rtuple$ we obtain isomorphisms $D_i(G_\rtuple) \cong 
N/N_\rtuple$ and hence $D_i(G_\rtuple)_\ab \cong N/[N,N]N_\rtuple$. As there are only finitely many possibilities for~$N$, it thus suffices to show that the set
$$\bigl\{r\in(N^-)^n \bigm| N/[N,N]N_\rtuple \hbox{ is infinite} \bigr\}$$
is closed of measure~$0$. But that is the content of Lemma \ref{DiabInfClosMeas0}.
\end{Proof}


\medskip
\begin{Proofof}{Theorem \ref{StrongThm}}
We claim that 
\UseTheoremCounterForNextEquation
\begin{equation}\label{StrongThmIntersect}
\Sch^\strong\ =\ \bigcap_{i\ge1}\; \bigl\{[G]\in\Sch \bigm| D_i(G)_\ab \hbox{ is finite} \bigr\}.\end{equation}
First observe that, since $G$ is finitely generated, every $D_i(G)$ is an open subgroup, proving the inclusion ``$\subset$''. 
To show the reverse inclusion consider any element $[G]$ of the right hand side. For any open subgroup $H<G$ there then exists an integer $i\ge1$ such that $D_i(G)\subset H$. Since $D_i(G)_\ab$ is finite, its commutator subgroup $[D_i(G),D_i(G)]$ is then an open subgroup of $D_i(G)$ and hence of~$H$. As it is contained in $[H,H]$, it follows that $[H,H]$ is an open subgroup as well. Thus $H_\ab=H/[H,H]$ is finite, and so $[G] \in \Sch^\strong$, as desired.

By taking complements, Proposition \ref{DiableI} implies that each term of the intersection in~\eqref{StrongThmIntersect} is an open subset of measure~$1$. Thus the intersection is a countable intersection of open subsets and itself has measure~$1$, as desired.
\end{Proofof}

%

\section{$p$-Tower groups}
\label{pTG}

Consider an imaginary quadratic field $K$ and let $K_p$ denote its maximal unramified pro-$p$-extension. Then $K_p$ is Galois over~$\BQ$ and its Galois group is the semidirect product of the pro-$p$-group $G_K := \Gal(K_p/K)$ with a group of order~$2$ generated by complex conjugation. Thus $G_K$ is a $\sigma$-pro-$p$-group. By Koch-Venkov \cite[\S1]{KochVenkov1974} and Proposition \ref{SchurSigmaProp} this is a strong Schur $\sigma$-group. Following McLeman \cite{McLeman2008} we call it the \emph{$p$-tower group} associated to~$K$.

\medskip
By abelian class field theory the maximal abelian quotient $G_{K,\ab}$ of $G_K$ is naturally isomorphic to the $p$-primary part of the narrow ideal class group of~$\CO_K$. Cohen and Lenstra \cite{CohenLenstra1983} described a probabilistic explanation for how often this group is isomorphic to a given finite abelian $p$-group and provided evidence for it. Our construction of the probability space $\Sch$ is intended to provide an analogue of this for the isomorphism class of~$G_K$.

\medskip
To explain this heuristic consider the set of imaginary quadratic fields up to isomorphism. 
(Perhaps one can also consider a natural infinite subset thereof, such as all those with discriminant in some congruence class, or which behave in a certain way vis-a-vis a fixed number field.)
Enumerate this set in some natural way, such as by increasing absolute value of the discriminant. The heuristic then says that for any suitable subset $P\subset\Sch$, the proportion of $K$ with $[G_K]\in P$ among the first $N$ fields should tend to $\mu_\infty(P)$ for $N\to\infty$. 

\medskip
As with applications of the Cebotarev density theorem to $\ell$-adic representations, without additional considerations
one can really only expect this to hold when the subset $P$ is open and closed. In particular, Proposition \ref{UiGMeasure} suggests that for any weak Schur $\sigma$-group $G$ and any integer $i\ge2$, the proportion of fields $K$ with $G_K/D_i(G_K)\cong G/D_i(G)$ should tend to the positive number
\UseTheoremCounterForNextEquation
\begin{equation}\label{UiGMeasureRepeat}
\frac{C_\infty}{C_{n-m}}\cdot \frac{1}{\bigl|\Aut_\sigma(G/D_i(G))\bigr|},
\end{equation}
where $n$ and $m$ depend in a specific way on $G$ and~$i$. Similar results can be deduced from this for other properties of $G_K$ that depend only on some finite quotient,
using Propositions \ref{UDGMeasureFin} and \ref{UDGMeasureFinVar}
In particular, Proposition \ref{FinSchurMeas} suggests that for any finite Schur $\sigma$-group~$G$, the proportion of fields $K$ with $G_K\cong G$ should tend to the positive number 
\UseTheoremCounterForNextEquation
\begin{equation}\label{FinSchurMeasRepeat}
\frac{C_\infty}{\bigl|\Aut_\sigma(G)\bigr|}.
\end{equation}

Note that by Corollary \ref{GAbMeasure} our measure gives the same values for the probability that $G_{K,\ab}$ is isomorphic to a given finite abelian $p$-group as in the original Cohen-Lenstra heuristic (see \cite{CohenLenstra1983} and \cite[Conj.\,5.10.1\,(3)]{Cohen1993}).

\medskip
One should be aware that for a subset $P$ of measure zero, the heuristic only says that the proportion of fields with $[G_K]\in P$ should tend to~$0$, but cannot rule out that this happens infinitely often. 

Also, the heuristic should not be expected to hold for all measurable subsets~$P$. This is related to the fact that, as there are only countably many isomorphism classes of imaginary quadratic fields, there are only countably many corresponding points $[G_K] \in\Sch$, while in all likelihood the measure space $\Sch$ has the cardinality of the continuum. 
Specifically, in \cite[Prop.\,4.7]{Pink2025} we will prove that $\mu_\infty(\{[G]\})=0$ for any infinite weak Schur $\sigma$-group~$G$. Thus the set of all $[G_K]$ with $G_K$ infinite also has measure~$0$, although $G_K$ should be infinite for a positive proportion of fields $K$ by results of Golod-Shafarevich type discussed below.

\medskip
Although $G_K$ is always a \emph{strong} Schur $\sigma$-group, the natural result of our construction in terms of generators and relations is a probability measure on the larger set $\Sch$ of isomorphism classes of all \emph{weak} Schur $\sigma$-groups. The heuristic can therefore only be plausible if the complement $\Sch\setminus\Sch^\strong$ has measure~$0$. That we can indeed prove this in Theorem \ref{StrongThm} indicates that our heuristic may be reasonable.

\medskip
A further fundamental question is under which conditions $G_K$ is finite, respectively infinite. After Golod and Shafarevich \cite{GolodShafarevich1964} proved that $G_K$ can be infinite, there has been much work on this, with \cite{Koch1969}, \cite{KochVenkov1974}, \cite{KochVenkov1975}, \cite{McLeman2008}, \cite{AhlqvistCarlson2025} among others. 
Many of these results can be phrased for abstract Schur $\sigma$-groups.

So let $G$ be a strong Schur $\sigma$-group with $d_G$ generators. According to current knowledge $G$ is finite if $d_G\le1$, infinite if $d_G\ge3$, and in the case $d_G=2$ it is infinite unless its Zassenhaus type is (3,3) or (3,5) or (3,7) (see McLeman \cite[\S2]{McLeman2008}). 
Also, finite groups of type (3,3) occur as $p$-tower groups by Scholz-Taussky \cite{ScholzTaussky1934} and Bartholdi-Bush \cite[Cor.\,3.3]{BartholdiBush2007}, and finite groups of type (3,5) and (3,7) occur by Ahlqvist-Carlson \cite[Thm.\,5.7]{AhlqvistCarlson2025}.
In particular there exist finite Schur $\sigma$-groups of these types.

Somewhat less seems to be known about infinite groups of these types as $p$-tower groups.
McLeman \cite[Conj.\,2.9]{McLeman2008} conjectured that any $p$-tower group of Zassenhaus type (3,3) should be finite. In view of our heuristic this suggests that most Schur $\sigma$-groups of Zassenhaus type (3,3) should be finite. That is a purely group theoretical statement and might therefore be easier to decide.
In fact, for $p>3$ it is proved in a separate paper by the first author \cite[Thm.\,8.5]{Pink2025}:

\begin{Thm}\label{Schur33Thm}
For $p>3$ the set of $[G]\in\Sch$ such that $G$ is infinite of Zassenhaus type (3,3) has measure~$0$.
\end{Thm}

We have no opinion about whether the same statement should hold for $p=3$. 




\end{document}